\documentclass{amsart}
\usepackage{times}
\usepackage{amsmath}
\usepackage{amsfonts}
\usepackage{latexsym}
\usepackage[T1]{fontenc}
\usepackage[latin1]{inputenc}
\usepackage{hyperref}
\usepackage[arrow,curve,matrix]{xy}
\usepackage{epsfig}
\DeclareMathOperator{\Alb}{Alb}
\DeclareMathOperator{\Aut}{Aut}
\DeclareMathOperator{\Chow}{Chow}
\DeclareMathOperator{\Hom}{Hom}

\DeclareMathOperator{\RatCurves}{RatCurves}
\DeclareMathOperator{\Sing}{Sing}

\newcommand{\branch}{{\mathfrak B}}

\newcommand{\codim}{{\rm codim}}
\newcommand{\Pic}{{\rm Pic}}

\newcommand{\ramification}{{\mathfrak R}}

\newtheorem{lemma1}[equation]{}

\theoremstyle{plain}
\newtheorem{thm}{Theorem}[section]
\newtheorem{defn}[thm]{Definition}

\numberwithin{equation}{thm}
\numberwithin{figure}{section}
\theoremstyle{plain}
\newtheorem{cor}[thm]{Corollary}

\newtheorem{conj}[thm]{Conjecture}
\newtheorem{lem}[thm]{Lemma}

\newtheorem{ex}[thm]{Example}

\theoremstyle{plain}
\newtheorem{prop}[thm]{Proposition}
\newtheorem{proclaim-special}[thm]{\specialthmname}

\theoremstyle{remark}
\newtheorem{fact}[thm]{Fact}
\newtheorem{rem}[thm]{Remark}
\newtheorem{subrem}[equation]{Remark}
\theoremstyle{remark}
\newtheorem*{claim*}{Claim}
\newtheorem{notation}[thm]{Notation}
\newtheorem{assumption}[thm]{Assumption}
\newtheorem{assnot}[thm]{Assumption / Notation}
\newtheorem{problem}[thm]{Problem}

\def\factor#1.#2.{\left. \raise 2pt\hbox{$#1$} \right/\hskip -2pt\raise
  -2pt\hbox{$#2$}}

\newcounter{thisthm}
\newcommand{\ilabel}[1]{\newcounter{#1}\setcounter{thisthm}{\value{thm}}\setcounter{#1}{\value{enumi}}}
\newcommand{\iref}[1]{(\thesection.\the\value{thisthm}.\the\value{#1})}

\newenvironment{enumerate-c}{
  \begin{enumerate}
    \setcounter{enumi}{\value{equation}}}
  {\setcounter{equation}{\value{enumi}}\end{enumerate}}

\newcommand\sE{{\mathcal E}}
\newcommand\sF{{\mathcal F}}
\newcommand\sG{{\mathcal G}}
\newcommand\sH{{\mathcal H}}
\newcommand\sI{{\mathcal I}}

\newcommand\sL{{\mathcal L}}
\newcommand\sO{{\mathcal O}}

\newcommand\bR{{\mathbb R}}

\newcommand\bQ{{\mathbb Q}}

\newcommand\tor{{\rm tor}}

\newcommand\Supp{{\rm Supp}}
\DeclareMathOperator{\rank}{rank}

\newcommand\bP{{\mathbb P}}

\def\factor#1.#2.{\left. \raise 2pt\hbox{$#1$} \right/\hskip -2pt\raise
  -2pt\hbox{$#2$}}

\newlength{\swidth}
\newenvironment{enumerate-p}{
  \begin{enumerate}}
  {\setcounter{equation}{\value{enumi}}\end{enumerate}}

\title{Endomorphisms of projective varieties}

\author{Marian Aprodu, Stefan Kebekus and Thomas Peternell}

\date{\today}

\thanks{Marian Aprodu was supported by a Humboldt Research Fellowship
  and a Humboldt Return Fellowship. He expresses his special thanks to
  the Mathematical Institute of Bayreuth University for hospitality
  during the first stage of this work. Stefan Kebekus and Thomas
  Peternell were supported by the DFG-Schwerpunkt ``Globale Methoden
  in der komplexen Geometrie'' and the DFG-Forschergruppe
  ``Classification of Algebraic Surfaces and Compact Complex
  Manifolds''. A part of this paper was worked out while Stefan
  Kebekus visited the Korea Institute for Advanced Study.  He would
  like to thank Jun-Muk Hwang for the invitation.}

\address{Marian Aprodu, Institute of Mathematics ``Simion Stoilow'' of
  the Romanian Academy, P.O. BOX 1-764, 014700 Bucharest, Romania and
  \c Scoala Normal\u a Superioar\u a, Calea Grivitei 21, 010702
  Bucharest, Romania}
\email{\href{mailto:Marian.Aprodu@imar.ro}{Marian.Aprodu@imar.ro}}
\urladdr{\href{http://www.imar.ro/~aprodu}{http://www.imar.ro/~aprodu}}

\address{Stefan Kebekus, Mathematisches Institut, Universit\"at zu
  K\"oln, Weyertal 86--90, 50931 K\"oln, Germany}
\email{\href{mailto:stefan.kebekus@math.uni-koeln.de}
{stefan.kebekus@math.uni-koeln.de}}
\urladdr{\href{http://www.mi.uni-koeln.de/~kebekus}
{http://www.mi.uni-koeln.de/~kebekus}}

\address{Thomas Peternell, Mathematisches Institut, Universit\"at 
Bayreuth, 95440 Bayreuth, Germany}
\email{\href{mailto:thomas.peternell@uni-bayreuth.de}
{thomas.peternell@uni-bayreuth.de}}

\begin{document}

\maketitle

\begin{abstract}
  We study complex projective manifolds $X$ that admit surjective
  endomorphisms $f:X\to X$ of degree at least two. In case $f$ is
  \'etale, we prove structure theorems that describe $X$. In
  particular, a rather detailed description is given if $X$ is a
  uniruled threefold. As to the ramified case, we first prove a
  general theorem stating that the vector bundle associated to a
  Galois covering of projective manifolds is ample (resp.\ nef) under
  very mild conditions.  This is applied to the study of ramified
  endomorphisms of Fano manifolds with $b_2 = 1$. It is conjectured
  that $\bP_n$ is the only Fano manifold admitting admitting an
  endomorphism of degree $d \geq 2$, and we prove that in several
  cases.

  A part of the argumentation is based on a new characterization of
  $\bP_n$ as the only manifold that admits an ample subsheaf in its
  tangent bundle.
\end{abstract}

\setcounter{tocdepth}{1}
\tableofcontents

\section{Introduction}

A classical question in complex geometry asks for a description of
projective manifolds $X$ that admit surjective endomorphisms $f:X\to
X$ of degree at least two ---we refer to \cite{Fakhruddin03} for
questions and results relating endomorphisms of algebraic varieties
with some general conjectures in number theory. A straightforward
argument proves that $X$ cannot be of general type and that $f$ is
necessarily finite, see \cite[Prop.~2]{Beauville01} and
\cite[Lem.~2.3]{Fu02}. While a complete classification of all possible
$X$ remains an extremely hard open problem, the present paper presents
a number of results. Two different approaches are highlighted in the
text, according to whether the given endomorphism is \'etale or not.

\subsection{\'Etale endomorphisms and the minimal model program}

The first approach, started in Section~\ref{sec:2}, is based on the
Minimal Model Program and deals mainly with \'etale endomorphisms.
This approach was also used by Fujimoto in \cite{Fu02}. One can check
without much difficulty that the direct image map $f_*: H_2(X,\mathbb
R) \to H_2(X,\mathbb R)$ permutes the geometric extremal rays of the
Mori cone, and the same is true for extremal rays if $f$ is \'etale,
cf.~Proposition~\ref{prop:2extr-2} and \cite[Prop.~4.2]{Fu02}.  The
targets of the associated contractions are related by finite
morphisms, cf.~Corollaries~\ref{cor:2contr}--\ref{cor:2simplyconnE-b},
Proposition \ref{prop:etalebetweennonisom} and \cite[Prop.~4.4]{Fu02}.
In Section~\ref{sec:3}, we prove any endomorphism of $X$ is be \'etale
if $K_X$ is pseudo-effective. This generalizes a previous result by
Iitaka \cite[Thm.~11.7]{Iitaka82}; see also \cite{Beauville01}. In
Section~\ref{sec:5}, we employ the results of Sections~\ref{sec:2} and
\ref{sec:3} to study \'etale endomorphisms in greater detail. We
observe that the existence of an \'etale endomorphism $f$ has some
implications for invariants of $X$.  For example, the top
self-intersection of the canonical bundle, the Euler characteristic,
and the top Chern class of the manifold must all vanish,
cf.~Lemma~\ref{lem:invariants}.  As a general result, we prove in
Proposition~\ref{prop:4deform} that all deformations of \'etale
endomorphisms come from automorphisms of $X$.  Furthermore, we study
the Minimal Model Program in detail if $X$ has dimension three.
Since the case on non-negative Kodaira dimension was carefully treated
in \cite{Fu02}, we focus on the case of negative Kodaira dimension.

\subsection{The vector bundle associated with an endomorphism}

The second approach is inspired by Lazarsfeld's work \cite{Laz80}, see
also \cite[6.3.D]{Laz-Positivity2}. The idea is to study a ramified
finite covering $f$ of degree $d$ through the properties of the
canonically associated vector bundle $\mathcal{E}_f$. Notably, the
bundle $\mathcal{E}_f$ tends to inherit positivity properties from the
ramification divisor. This works particularly well for Galois covers,
as shown in the following Theorem.

\begin{thm}[\protect{cf.~Theorem~\ref{thm:ampleness1}}]
  Let $f: X \to Y$ be a Galois covering of projective manifolds which
  does not factor through an \'etale covering.  Assume that all
  irreducible components of the ramification divisor are ample. Then
  $\sE_f$ is ample. \qed
\end{thm} 

The condition that $f$ does not factor through an \'etale covering is
automatically satisfied, e.g., for Fano manifolds $X$ with $b_2(X) =
1$. Notice that the theorem is false without the Galois assumption.

The approach via the bundle $\sE_f$ is applied in the last sections of
the paper, where we investigate manifolds of negative Kodaira
dimension admitting a ramified endomorphism. It is generally believed
that the projective space is the only Fano manifold with Picard number
one for which such non-trivial endomorphisms exist.

\begin{conj}\label{conj:12}
  Let $f: X \to X$ be an endomorphism of a Fano manifold $X$ with
  $\rho(X) = 1$. If $\deg f > 1, $ then $X \simeq \bP_n$.
\end{conj} 

At present, Conjecture~\ref{conj:12} known to be true in the following
special cases: surfaces, threefolds \cite{ARV99, Schuhmann99,
  Hwang-Mok03}, rational homogeneous manifolds \cite{PS89, HM99}, or
toric varieties \cite{OW02}, varieties containing a rational curve
with trivial normal bundle \cite[Cor.~3]{Hwang-Mok03}. In
Section~\ref{sec:branched}, we enlarge the list. In particular, we
prove the following results.
\begin{thm}[Indirect evidence for Conjecture~\ref{conj:12}]
  Let $X$ be a Fano manifold with $\rho(X) = 1$ and $f: X \to X$ an
  endomorphism.  If one of the following conditions holds, then $\deg
  f = 1$.
  \begin{itemize}
  \item $X$ has index $\leq 2$ and additionally there exists a line in
    $X$ which is not contained in the branch locus of $f$ ---see
    Theorem~\ref{thm:smallIndex}
  \item $X$ satisfies the Cartan-Fubini condition, $X$ is almost
    homogeneous and $h^0(X,\, T_X) > \dim X$ ---see
    Theorem~\ref{thm:CF}
  \item $X$ satisfies the Cartan-Fubini condition, $X$ is almost
    homogeneous and either branch or the ramification divisor of $f$
    meets the open orbit of $\Aut^0(X)$ ---see Theorem~\ref{thm:CF}
  \item $X$ is a del Pezzo manifold of degree $5$ and $\rho(X) = 1$
    ---see Theorem~\ref{thm:delPezzo}. \qed
  \end{itemize}  
\end{thm}

\begin{thm}[Direct evidence for Conjecture~\ref{conj:12}]
  Let $X$ be a Fano manifold with $\rho(X) = 1$ and $f: X \to X$ an
  endomorphism of degree $\deg f \geq 2$. Then $X \simeq \bP_n$ if one
  of the following conditions hold.
  \begin{itemize} 
  \item $\sE_{f_k}$ is ample for some iterate $f_k$ of $f$ and
    $h^0(X,\,f^*_k(T_X)) > h^0(X,\, T_X)$ ---see
    Corollary~\ref{cor:5critA}
  \item $f$ is Galois and $h^0(X,\,f^*_k(T_X)) > h^0(X,\,T_X)$ for
    some $k$ ---see Corollary~\ref{cor:611}
  \item $X$ is almost homogeneous and $h^0(X,\, T_X) > \dim X$ and
    $\sE_{f_k}$ is ample for sufficiently large $k$ ---see
    Corollary~\ref{cor:5critA} and Remark~\ref{cor:asmallx}. \qed
  \end{itemize} 
\end{thm} 

In order to check that $X$ is isomorphic to a projective space, we
prove the following partial generalization of a theorem of Andreatta
and Wisniewski, \cite{AW01}.
\begin{thm}[\protect{cf.~Theorem~\ref{thm:charactPn}}] 
  Let $X$ be a projective manifold with $\rho (X) = 1$. Let $\sF
  \subset T_X$ be a coherent subsheaf of positive rank. If $\sF$ is
  ample, then $X \simeq \bP_n$. \qed
\end{thm}

In fact, we prove a much stronger theorem, assuming the ampleness of $\sF$ only on certain
rational curves.

\section{Notation and general facts}
\label{sec:1e}

We collect some general facts on surjective endomorphisms.  Unless
otherwise noted, we fix the following assumptions and notation
throughout the present work.

\begin{assnot}\label{ass:endomorphism}
  Let $X$ be a projective manifold and $f: X \to X$ a surjective
  endomorphism. Let $d$ be the degree of $f$. The ramification divisor
  upstairs is denoted by $\ramification$, so that
  \begin{equation}\label{eq:adjunction_for_f}
    K_X = f^*(K_X) + \ramification.    
  \end{equation}
  The branch divisor downstairs is called $\branch$. It is defined as
  the cycle-theoretic image $\branch = f_*(\ramification)$.
\end{assnot}

We briefly recall two lemmas that hold on every compact manifold and
do not require any projectivity assumption.

\begin{lem}[\protect{\cite[Lem.~1]{Beauville01}, \cite[Lem.~2.3]{Fu02}}]
\label{lem:Beauville1}
  The linear maps
  $$
  f^*: H^*(X,\bQ) \to H^*(X,\bQ) \text{\quad and \quad}
  f_*:H^*(X,\bQ) \to H^*(X,\bQ)
  $$
  are isomorphisms. More precisely, $f_*f^* = d \cdot {\rm id}$,
  $f^*f_* = d \cdot {\rm id}$, and therefore $(f^*)^{-1} = \frac{1}{d}
  f_*$. In particular, $f$ is finite.  \qed
\end{lem}

\begin{lem}[\protect{\cite[Prop.~2]{Beauville01}}]
  Under the Assumption~\ref{ass:endomorphism}, the manifold $X$ is not
  of general type, i.e., $\kappa (X) < \dim X$. \qed
\end{lem}

\part*{\'Etale endomorphisms and the minimal model program}

\section{Extremal contractions of manifolds with endomorphisms}
\label{sec:2}

\subsection{Extremal contractions}
\label{sec:extrcontrwithend}

Recall that $N^1(X) \subset H^2(X,\bR)$ is the subspace generated by
the classes of irreducible hypersurfaces and that ${\overline {NE}}(X)
\subset N_1(X,\bR) \subset H_2(X,\bR)$ is the closed cone generated by
classes of irreducible curves. Equivalently,
$$
{\overline {NE}}(X) = \bigl\{ \alpha \in N_1(X,\bR) \,|\, H.\alpha
\geq 0 \text{ for all ample } H \in \Pic(X) \bigr\}.
$$

\begin{notation}
  A half-ray $R \subset {\overline {NE}}(X)$ is called \emph{extremal}
  if it is geometrically extremal and if $K_X \cdot R < 0$. We say
  that $\alpha \in {\overline {NE}}(X)$ is extremal if the ray $R =
  \bR_+ \alpha$ is extremal.
\end{notation}

We recall a proposition that has already been shown by Fujimoto. For
the convenience of the reader, a short argument is included.

\begin{prop}[\protect{cf.~\cite[Prop.~4.2]{Fu02}}]\label{prop:2extr-2}
  Under the Assumption~\ref{ass:endomorphism}, the linear isomorphism
  of $\bR$-vector spaces, $f_*: H_2(X,\bR) \to H_2(X,\bR)$, restricts
  to a bijective map
  $$
  f_*: {\overline {NE}}(X) \to {\overline {NE}}(X).    
  $$
  In particular, $f_*$ defines a bijection on the set of geometrically
  extremal rays. 
  \begin{enumerate}
  \item\ilabel{il:a} If $R \subset {\overline {NE}}(X)$ is an extremal
    ray such that the exceptional locus $E_R$ of the associated
    contraction is not contained in the ramification locus
    $\ramification$, then $f_*(R)$ is again extremal.
  \item\ilabel{il:b} If $f$ is \'etale, then $f_*$ defines a bijection
    on the set of extremal rays.
  \end{enumerate}
\end{prop}
\begin{proof}
  To show that $f_*$ defines a bijection on the set of geometrically
  extremal rays, it suffices to note that both the cycle-theoretic
  image and preimage of an irreducible, effective curve under the
  finite morphism $f$ is effective.
  
  For Statement~\iref{il:a}, let $R \subset {\overline {NE}}(X)$ be an
  extremal ray with associated exceptional set $E_R \subset X$. If
  $E_R \not \subset \ramification$, let $C \subset E_R$ be a curve
  with $[C] \in R$ and $C \not \subset \ramification$. Then
  $$
  f_*([C]).K_X = [C].f^*(K_X) = \underbrace{K_X\cdot C}_{< 0} -
  \underbrace{\ramification \cdot C}_{\geq 0} < 0.
  $$
  It follows that $f_*(R)=\mathbb R^+ \cdot f_*\bigl([C]\bigr)$ is
  extremal. This shows the second statement.
  
  For~\iref{il:b}, assume that $f$ is \'etale. We will only need to
  show that the pull-back of an extremal curve is extremal. To that
  end, let $R$ be an extremal ray. It is them immediately clear that
  $f^*(R)$ is geometrically extremal. To show that it is extremal, let
  $C \subset X$ be an irreducible rational curve with $[C] \in R$.
  Since
  $$
  K_X \cdot f^*([C]) = f^*(K_X) \cdot f^*([C]) = d \cdot (K_X \cdot
  [C]) < 0,
  $$
  the pull-back ray $f^*(R) = \mathbb R^+ \cdot f^*\bigl([C]\bigr)$
  is thus indeed extremal.
\end{proof}

\begin{rem}
  Under the Assumptions~\ref{ass:endomorphism} suppose that $X$ has
  only finitely many extremal rays. Then there exists a number $k$,
  with the following property: if $f_k$ is the $k^{\rm th}$ iteration
  of $f$ and if $R$ is any extremal ray, then $(f_k)_*(R) = R$.
  Namely, if $M$ is the finite set of extremal rays, then $(f_*)^k =
  (f_k)_*: M \to M$ is bijective. Hence there exists a number $k$ such
  that $f_*^k$ is the identity.
\end{rem}

\begin{cor}\label{cor:2contr}
  Let $R \subset \overline {NE}(X)$ be an extremal ray. If the
  associated contraction is birational, assume that the exceptional
  set $E_R$ of the associated contraction is not contained in
  $\ramification$.  Then there exists a commutative diagram as
  follows.
  $$
  \xymatrix{
    X \ar[d]_{f, \text{ finite}} \ar[rrr]^{\phi, \text{ contraction of $R$}}
    & & & Y  \ar[d]^{g, \text{ finite}} \\
    X \ar[rrr]_{\psi, \text{ contraction of $f_*(R)$}} & & & Z
  }
  $$
  In particular, if $E_{f_*(R)}$ is the exceptional set of $\psi$,
  then $f(E_R) = E_{f_*(R)}$ and $E_R = f^{-1}\bigl(E_{f_*(R)}\bigr)$.
\end{cor}
\begin{subrem}
  The $\phi$-exceptional set $E_R$ is defined as the set where $\phi$
  is not locally isomorphic. If $\phi$ is of fiber type, then $E_R =
  X$.
\end{subrem}
\begin{proof}
  Observe that if $\phi$ contracts a curve $C \subset X$, then $\psi$
  contracts $f(C) \subset X$, because the class $[f(C)]$ is contained
  in $f_*(R)$. Using Zariski's main theorem, this already shows the
  existence of $g$ and proves that $f(E_R) \subset E_{f_*(R)}$. Since
  $b_2(Y)=b_2(Z)=b_2(X)-1$, $g$ is necessarily finite.
  
  To show that $f(E_R) \supset E_{f_*(R)}$ and $E_R =
  f^{-1}\bigl(E_{f_*(R)}\bigr)$, let $C' \subset X$ be a curve which
  is contracted by $\psi$. Its class $[C']$ is then contained in
  $f_*(R)$, and $f^{-1}(C')$ is a union of curves whose individual
  classes are, by geometric extremality, contained in $R$. In
  particular, all components of $f^{-1}(C')$ are contained in
  $E_R$.
\end{proof}

\begin{cor}\label{cor:2simplyconnE-a}
  In the setup of Corollary~\ref{cor:2contr}, assume that $f: X \to X$
  is \'etale. If $\dim X \geq 5$, assume additionally that the
  contraction $\phi$ is divisorial.  Then the restriction $f|_{E_R}:
  E_R \to E_{f_*(R)}$ is \'etale of degree $d$. In particular,
  $E_{f_*(R)}$ is not simply connected.
\end{cor}
\begin{proof}
  If $\phi$ is of fiber type, there is nothing to show. We will thus
  assume that $\phi$ is birational. Observe that if the contraction
  $\phi$ is divisorial, the statement follows from
  Corollary~\ref{cor:2contr} and from the fact that the exceptional
  divisor of a divisorial contraction is irreducible. We consider the
  possibilities for $\dim X$.
  
  If $\dim X \leq 3$, it follows from the classification of extremal
  contractions in dimension 3, \cite[Thm.~3.3, Thm.~3.5]{Mori82}, that
  $\phi$ is divisorial. The claim is thus shown.
  
  If $\dim X = 4$, we are again finished if that show that $E_R$ is a
  divisor. We assume to the contrary and suppose $\dim E_R < 3$.  In
  this setup, a theorem of Kawamata \cite[Thm.~1.1]{Kaw89}, asserts
  that the exceptional loci of both $E_R$ and $E_{f_*(R)}$ are
  disjoint copies of $\bP_2$'s. In particular, they are simply
  connected. So, if $E_{f_*(R)}$ has $n$ connected components, then
  $E_R = f^{-1}\bigl(E_{f_*(R)}\bigr)$ will have $d\cdot n$ connected
  components. But the same argumentation applies to the \'etale
  morphism $f\circ f: X\to X$ and yields that $E_R$ has $d^2\cdot n$
  connected components. Again, we found a contradiction.
  
  Finally, assume that $\dim X \geq  5$. Then $\phi$ is divisorial by
  assumption.
\end{proof}

\begin{cor}\label{cor:2simplyconnE-b}
  In the setup of Corollary~\ref{cor:2simplyconnE-a}, if $\dim X = 2$,
  then $\phi$ is a ${\mathbb P}_1$-bundle. If $\dim X = 3$ and if
  $\phi$ is birational, then $g: Y \to Z$ is \'etale of degree $d$ and
  both $\phi$ and $\psi$ are blow-ups of smooth curves.
\end{cor}
\begin{proof}
  Recall that extremal loci of birational surface contractions are
  irreducible, simply connected divisors. This shows that $X$ is
  minimal, that $\phi$ is of fiber type and settles the case $\dim X =
  2$.
  
  Now assume that $\dim X = 3$ and that $\phi$ is birational. Observe
  that if $\phi(E_R)$ was a point, the classification
  \cite[Thm.~3.3]{Mori82} yields that $E_{f_*(R)}$ is isomorphic to
  $\bP_2$, $\bP_1 \times \bP_1$ or to the quadric cone.  But all these
  spaces are simply connected, a contradiction.  Consequence:
  $\phi(E_R)$ is not a point, and \cite[Thm.~3.3]{Mori82} asserts that
  $\phi$ is a blow-up.  The same holds for $\psi$, the contraction of
  $f_*(R)$. The fact that fibers of $\psi$ are 1-connected yields the
  \'etalit\'e of $g$.
\end{proof}

\subsection{Extremal contractions in the presence of \'etale morphisms between non-isomorphic varieties}

We remark that the results of
Proposition~\ref{prop:2extr-2}--Corollary~\ref{cor:2simplyconnE-b}
remain true for \'etale morphisms between possibly non-isomorphic
varieties, as long as their second Betti numbers agree.

\begin{prop}\label{prop:etalebetweennonisom}
  Let $g: X \to X'$ be a surjective \'etale morphism of degree $d$
  between projective manifolds that satisfy $b_2(X) = b_2(X')$.  Then
  the following holds.
  \begin{enumerate}
  \item The linear isomorphism of $\bR$-vector spaces, $g_*:
    H_2(X,\bR) \to H_2(X',\bR)$, restricts to a bijective map $ g_*:
    {\overline {NE}}(X) \to {\overline {NE}}(X') $ and defines a
    bijection on the set of extremal rays.
    
  \item If $R \subset \overline {NE}(X)$ is an extremal ray, then
    there exists a commutative diagram as follows
    $$
    \xymatrix{
      X \ar[d]_{g, \text{ \'etale}} \ar[rrr]^{\phi, \text{
          contraction of $R$}} & & & Y  \ar[d]^{h, \text{ finite}} \\
      X' \ar[rrr]_{\psi, \text{ contraction of $f_*(R)$}} & & & Z.  }
    $$
    If $E_R$ and $E_{g_*(R)}$ are the exceptional sets of $\phi$
    and $\psi$, respectively, then $g(E_R) = E_{g_*(R)}$ and $E_R =
    g^{-1}(E_{g_*(R)})$. If the contraction $\phi$ is divisorial, then
    $g|_{E_R}: E_R \to E_{g_*(R)}$ is \'etale of degree $d$ and
    $E_{g_*(R)}$ is thus not simply connected.
  
  \item If $\dim X = 2$, then $X$ is minimal.  If $\dim X = 3$ and
    $\phi$ is birational, then $h$ is \'etale of degree $d$ and both
    $\phi$ and $\psi$ are blow-ups of smooth curves.
  \end{enumerate}
\end{prop}
\begin{proof}
  Since $g_*\circ g^*: H^*(X') \to H^*(X)$ is multiplication with $d =
  \deg(g)$, we have that $g^*: H^*(X') \to H^*(X)$ is (piecewise)
  injective. Since $b_2(X) = b_2(X')$, we deduce as in
  Lemma~\ref{lem:Beauville1} that both $g^*: H^2(X') \to H^2(X)$ and
  $g_*: H^2(X) \to H^2(X')$ are isomorphic. The argumentation of
  Section~\ref{sec:extrcontrwithend} can then be applied verbatim.
\end{proof}

Proposition~\ref{prop:etalebetweennonisom} will later be used in the
following context. In the setup of Corollary~\ref{cor:2contr}, assume
that $\phi: X \to Y$ is the blow-up of the projective manifold $Y$
along a smooth curve $C$. Let $E = \phi^{-1}(C)$. It is then not
difficult to see that $\psi: X \to Z$ is then also a blow-up of a
manifold along a curve. Since $b_2(Y) = b_2(Z)$,
Proposition~\ref{prop:etalebetweennonisom} applies to contractions of
$Y$.

\section{Endomorphisms of non-uniruled varieties}
\label{sec:3}

In \cite[Thm.~11.7,~p.~337]{Iitaka82} it was shown that endomorphisms
of manifolds with $\kappa (X) \geq 0$ are necessarily \'etale. We
generalize this to non-uniruled varieties, at least when $X$ is
projective. First we state the following weaker result which also
holds for K\"ahler manifolds.

\begin{thm}\label{thm:3main}
  Let $X$ be a compact K\"ahler manifold and $f: X \to X$ a surjective
  endomorphism. If $K_X$ is pseudo-effective, i.e.~if its class is in
  the closure of the K\"ahler cone, then $f$ is \'etale.
\end{thm}

\begin{proof}
  We argue by contradiction: assume that $K_X$ is pseudo-effective and
  $f$ not \'etale, i.e.~assume that the ramification divisor of $f$ is
  not trivial: $\ramification \ne 0$. We fix a K\"ahler form $\omega$
  on $X$. It is an immediate consequence of Lemma~\ref{lem:Beauville1}
  that $f$ is finite. The standard adjunction formula for a branched
  morphism, $K_X = f^*(K_X) + \ramification$, then has two
  consequences:

  First, the canonical bundle is not numerically trivial, $K_X \not
  \equiv 0$. Since $K_X$ is assumed pseudo-effective, that means $K_X
  \cdot \omega^{n-1} > 0$ (in fact, the class of $K_X$ is represented by a non-zero
positive closed current $T$ and it is a standard fact that $T \cdot \omega^{n-1} > 0$
unless $T = 0$). 

  Secondly, if $f_m := f \circ \cdots \circ f$ is the $m^{\rm th}$
  iteration of $f$, the iterated adjunction formula reads
  $$
  K_X = f_m^*(K_X) + f_{m-1}^*(\ramification) + \ldots + f^*(\ramification) + \ramification.
  $$
  Intersecting with $\omega^{n-1}$, we obtain
  \begin{equation}\label{eq:3itint}
    K_X \cdot \omega^{n-1} = f_m^*(K_X) \cdot
    \omega^{n-1} + f_{m-1}^*(\ramification) \cdot \omega^{n-1} + \ldots + \ramification \cdot \omega^{n-1}.
  \end{equation}
  Observe that there exists a number $c > 0$ such that $L\cdot
  \omega^{n-1} > c$ for all non-trivial pseudo-effective line bundles
  $L \in \Pic(X)$; the number $c$ exists because the cohomology
  classes of pseudo-effective line bundles are exactly the integral
  points in the pseudo-effective cone. Since all of the $m+1$ summands
  in equation~\eqref{eq:3itint} are therefore larger than $c$, we have
  $$
  K_X \cdot \omega^{n-1} \geq (m+1) c
  $$
  for all positive integers $m$. This is absurd.
\end{proof}

\begin{cor}\label{cor:3branchedgiveuniruled}
  Let $X$ be a projective manifold, $f: X \to X$ surjective. If $X$ is
  not uniruled, then $f$ is \'etale.
\end{cor}

\begin{proof}
  Since $X$ is not uniruled, $K_X$ is pseudo-effective by
  \cite{BDPP04}.  Now apply Theorem~\ref{thm:3main}.
\end{proof}

\begin{rem}
  If $X$ is K\"ahler and $f: X \to X$ surjective but not \'etale,
  Theorem~\ref{thm:3main} asserts that $K_X$ cannot be
  pseudo-effective. It is, however, unknown whether this implies that
  $X$ is uniruled. For that reason we cannot state
  Corollary~\ref{cor:3branchedgiveuniruled} in the K\"ahler case
  although we strongly believe that it will be true.
  
  The proof of Theorem~\ref{thm:3main} shows a little more: if
  $\omega$ is any K\"ahler form and $\eta$ any positive closed
  $(n-1,n-1)$-form, we have $ K_X \cdot \omega^{n-1} \leq 0$ and $K_X
  \cdot \eta \leq 0 $. Since $K_X \not \equiv 0$, we must actually
  have strict inequality for some $\eta$. This could be useful in a
  further study of the K\"ahler case.
\end{rem}

\section{\'Etale endomorphisms}
\label{sec:5}

In this section we study \'etale endomorphisms more closely. In
Section~\ref{sec:5a}, we consider varieties of arbitrary dimension,
and study endomorphisms from a deformation-theoretic point of view.
In Section~\ref{sec:5balb} we study the interaction of the
endomorphism with the Albanese map.  Next we restrict to dimension 3
and apply the minimal model program to $X$.  Since the case of
non-negative Kodaira dimension was treated by Fujimoto, we restrict
ourselves to uniruled threefolds $X$.

Maintaining the Assumptions~\ref{ass:endomorphism}, we suppose
throughout this section that $f$ is \'etale. We set $n := \dim X$ and
note that a number of invariants vanish.

\begin{lem}\label{lem:invariants}
  In this setup, we have the following numerical data.
  \begin{enumerate}
  \item The class $c_1(K_X)^n \in H^{2n}(X)$ is zero.
  \item $\chi(\sO_X) = 0$.
  \item $c_n(X) = 0$.
  \end{enumerate}
\end{lem}
\begin{proof} 
  The first claim follows from $c_1(K_X)^n = f^*(c_1(K_X)^n) = d \cdot
  c_1(K_X)^n$. For the second claim, observe $\chi(\sO_X) = d \cdot
  \chi(\sO_X)$. The third results from $T_X = f^*(T_X)$.
\end{proof}

\subsection{Deformations of \'etale endomorphisms}\label{sec:5a}

This section is concerned with a study of deformations of $f$. We will
show that all deformations of the \'etale morphism $f$ come from
automorphisms of $X$.  This strengthens the results of \cite{HKP03,
  KP05} in our case.

\begin{thm}\label{prop:4deform}
  For an \'etale morphism $f: X \to X$, we have
  \begin{equation}\label{eq:haut}
    \Hom_f(X,X) \cong \factor \Aut^\circ(X) . \Aut^\circ(X) \cap
    \Aut(X/X).,
  \end{equation}
  where $\Aut^\circ(X)$ is the maximal connected subgroup of
  $\Aut(X)$ and $\Hom_f(X,X)$ is the connected component of
  $\Hom(X,X)$ that contains $f$. In particular, $\Hom_f(X,X)$ is
  irreducible, reduced and smooth.
  
  If $X$ is not uniruled and if $h^0(X, T_X) > 0, $ then there exists
  a finite \'etale cover $g: \tilde X \to X$ where $\tilde X$ is a
  product $\tilde X = A \times W$ such that $A$ is a torus and such
  that $H^0(A, T_A) = g^*H^0(X, T_X)$.
\end{thm}
\begin{proof}
  Consider the quasi-finite composition morphism of quasi-projective
  schemes
  \begin{equation}\label{eq:compmor}
  \begin{array}{rccc}
    f^\circ : & \Aut^\circ(X) & \to & \Hom_f(X,X) \\
    & g & \mapsto & f\circ g
  \end{array}
  \end{equation}
  If we identify the tangent spaces $T_{\Aut^\circ(X)} \cong H^0(X,
  T_X)$, and $T_{\Hom_f(X,X)} \cong H^0(X, f^*(T_X))$, then the
  tangent morphism of $f^\circ$ at $Id \in \Aut^\circ(X)$, is simply
  the tangent map of $f$, i.e.,
  $$
  Tf^\circ|_{Id} = H^0(Tf) : H^0(X, T_X) \to H^0(X, f^*(T_X)).
  $$
  Since $f$ is \'etale, $Tf^\circ|_{Id}$ is necessarily isomorphic.
  Using the group structure of $\Aut^\circ(X)$, the same argument
  yields that $Tf^\circ|_{g}$ is isomorphic for all $g \in
  \Aut^\circ(X)$. In particular, since $\Aut^\circ(X)$ is reduced and
  smooth, the image $U := f^\circ(\Aut^\circ(X))$ of the morphism is
  an open neighborhood of $f$ in the Hom-scheme, and $f^\circ$
  identifies $U$ with the right hand side of~\eqref{eq:haut}. The
  bijectivity of the tangent map yields that $\Hom_f(X,X)$ is reduced
  and smooth along $U$.
  
  It remains to show that $f^\circ$ is set-theoretically surjective.
  If not, let $f' \in \partial U \setminus U$ be a point in the
  boundary. 
  
  We claim that morphism $f': X \to X$ is then again surjective and
  \'etale. For surjectivity, observe that the morphisms $f$ and $f'$
  are homotopic, so that the associated pull-back morphisms on
  cohomology are equal. But a proper morphism surjective if and only
  if the pull-back of the orientation form is non-zero. For
  \'etalit\'e, consider the morphism
  $$
  \eta : \Hom_f(X,X) \to \Pic(X), \text{\quad where \quad} \eta(g) 
  := g^*(K_X)-K_X.
  $$
  Observe that the image of $\eta(g)$ is the trivial bundle iff $g$
  is \'etale, and that $\eta$ is constant on the open set $U \subset
  \Hom_f(X,X)$.
  
  Now $f'$ being surjective and \'etale, we can again consider the
  composition morphism $(f')^\circ$, defined in analogy
  with~\eqref{eq:compmor}.  Its image $(f')^\circ(\Aut^\circ(X))$ is
  again open and therefore necessarily intersects $U$. This is to say
  that there are automorphisms $g, g' \in \Aut^\circ(X)$ such that
  $$
  f\circ g = f' \circ g', \quad \text{that is,} \quad f' = f \circ
  (g\circ (g')^{-1}),
  $$
  a contradiction to $f' \not \in U$.
  
  For the second statement, observe that a vector field on $X$ cannot
  have a zero, because otherwise $X$ would be uniruled. Hence
  \cite{Lib78} gives the decomposition.
\end{proof}

\subsection{The Albanese map of a variety with \'etale endomorphisms}
\label{sec:5balb}

We maintain the assumption that $f$ is \'etale and study the Albanese
map, see also \cite{Fu02}.

\begin{prop}
  Let $\alpha: X \to \Alb(X)$ be the Albanese and $Y \subset \Alb(X)$
  its image.  Then $f$ induces a finite \'etale cover $Y \to Y$.
  
  There exists a morphism $h: Y \to W$ to a variety of general type
  which is a torus bundle with fiber $B$ and which is trivialized
  after finite \'etale cover of $W$. The map $W \to W$ induced by $f$
  is an automorphism.
\end{prop}

\begin{proof}
  The universal property of the Albanese implies that $f$ induces an
  \'etale morphism $g: \Alb(X) \to \Alb(X)$ mapping $Y$ to $Y$.
  
  If $Y = \Alb(X)$, let $W$ be a point. If $Y \ne \Alb(X)$, we obtain
  a map $h: Y \to W$ where $W$ is of general type and the fibers are
  translates of subtori of $\Alb(X)$. Consider the induced map $g': W
  \to W$, whose existence is guaranteed by a classical result of
  Ueno, see e.g. \cite[Thm.~3.7]{Mori87}.  Since $W$ is of general
  type, i.e., since any desingularization is of general type, $g'$ is
  easily seen to be an automorphism --- adapt the proof
  \cite[Prop.~2]{Beauville01}. The fact that $h$ can be trivialized by
  a finite \'etale cover of $W$ is again Ueno's theorem cited above.
\end{proof}

\subsection{Extremal contractions of threefolds with \'etale endomorphisms}

As in \cite[Sect.~4]{Fu02}, we will now investigate \'etale morphisms
of threefolds more closely. It will turn out without much work that
any extremal contraction is the blow-up of a curve. More precisely,
the following strengthening of
Proposition~\ref{prop:etalebetweennonisom} holds true.

\begin{prop}\label{prop:divisorialcontr}
  Maintaining the Assumptions~\ref{ass:endomorphism}, suppose that $X$
  is a 3-fold and that $f$ is \'etale of degree $d$. Then there exists
  a commutative diagram
  \begin{equation}\label{eq:MMP}
    \xymatrix{
      X \ar[r]^{\phi_0} \ar[d]_{f} & Y_1 \ar[r]^{\phi_1} \ar[d]^{h_1} & 
      Y_2 \ar[r]^{\phi_1} \ar[d]^{h_2} & \cdots \ar[r]^{\phi_{m-1}} & 
      Y_m \ar[d]^{h_m} \\
      X \ar[r]_{\psi_0} & Z_1 \ar[r]_{\psi_1} & Z_2 \ar[r]_{\psi_2} & 
      \cdots \ar[r]_{\psi_{m-1}} & Z_m }
  \end{equation}
  where the $h_i$ are \'etale of degree $d$, the $\phi_i$ and $\psi_i$
  are extremal contractions and the following holds.
  \begin{enumerate-p}
  \setcounter{enumi}{\value{equation}}
  \item All $Y_i$ and $Z_i$ are smooth, $\phi$ and $\psi$ are blow-ups
    along smooth elliptic curves on which the canonical bundle is
    numerically trivial.
  \item Either $K_{Y_m}$ and $K_{Z_m}$ are nef or $Y_m$ and $Z_m$
    admit only contractions of fiber type.
  \end{enumerate-p}
\end{prop}

\begin{proof}
  Applying Corollary~\ref{cor:2contr} and
  Proposition~\ref{prop:etalebetweennonisom} inductively, we obtain an
  infinite diagram, with \eqref{eq:MMP} as the first two rows.
  \begin{equation}\label{eq:arrofofcontrs}
    \xymatrix{
      X \ar[r]^{\phi_0} \ar[d]_{f} & Y_1 \ar[r]^{\phi_1} \ar[d]^{h_{1,0}} &
      Y_2 \ar[r]^{\phi_1} \ar[d]^{h_{2,0}} & \cdots \ar[r]^{\phi_{m-1}} &
      Y_m \ar[d]^{h_{m,0}} \\
      X \ar[r]_{\psi_{0,1}} \ar[d]_{f} & Z_{1,1} \ar[r]_{\psi_{1,1}}
      \ar[d]^{h_{1,1}} & Z_{2,1} \ar[r]_{\psi_{2,1}} \ar[d]^{h_{2,1}} &
      \cdots \ar[r]_{\psi_{m-1,1}} & Z_{m,1} \ar[d]^{h_{m,1}} \\
      X \ar[r]_{\psi_{0,2}} \ar[d]_{f} & Z_{1,2} \ar[r]_{\psi_{1,2}}
      \ar[d]^{h_{1,2}} & Z_{2,2} \ar[r]_{\psi_{2,2}} \ar[d]^{h_{2,2}} &
      \cdots \ar[r]_{\psi_{m-1,2}} & Z_{m,2} \ar[d]^{h_{m,2}} \\
      \vdots &\vdots &\vdots & &\vdots  }
  \end{equation}
  The \'etalit\'e of the $h_{i,j}$, the smoothness of the $Z_{i,j}$
  and fact that the $\psi_{i,j}$ are blow-ups of smooth curves $D_i
  \subset Y_i$ and $C_{i,j} \subset Z_{i,j}$ are immediate from
  Proposition~\ref{prop:etalebetweennonisom}.  Statement (2) of
  Proposition~\ref{prop:divisorialcontr} is a standard result of
  minimal model theory.
    
  It remains to show that the curves $D_i$ and $C_{i,j}$ are elliptic.
  To this end, observe that the $h_{i,j}$ restrict to \'etale morphisms
  of the $C_{i,j}$, of degree $d$. Thus, for any number $i$ and $j$,
  $$
  \deg K_{D_i} = d^j \cdot \deg K_{C_{i,j}} \text{ and } \deg
  K_{Y_i}|_{D_i} = d^j \cdot \deg K_{Z_{i,j}}|_{C_{i,j}} .
  $$
  This is possible only if $\deg K_{D_i} = \deg K_{Y_i}|_{D_i} =
  0$. The same argument holds for any $C_{i,j}$.
\end{proof}

\begin{notation}
  In the setup of Proposition~\ref{prop:divisorialcontr}, we call
  $h_m: Y_m \to Z_m$ a \emph{minimal model} of $f$.
\end{notation}

\subsection{The minimal model program if $\kappa(X)=-\infty$}

Since the case $\kappa (X) \geq 0$ was studied in \cite{Fu02} in great
detail, we are mainly interested in the case where $\kappa(X) = -
\infty$. In this setup the manifolds $Y_m$ and $Z_m$ of
Proposition~\ref{prop:divisorialcontr} allow extremal contractions of
fiber type. To fix notation, we summarize the obvious properties in
the following Lemma.

\begin{lem}\label{lemnot:fibercontr}
  Maintaining the assumptions of
  Proposition~\ref{prop:divisorialcontr}, suppose that $\kappa(X) = -
  \infty$. Then the commutative diagram~\eqref{eq:MMP} extends as
  follows.
  $$
  \xymatrix{
    X \ar[r]^{\phi_0} \ar[d]_{f} & Y_1 \ar[r]^{\phi_1} \ar[d]^{h_1} &
    Y_2 \ar[r]^{\phi_1} \ar[d]^{h_2} & \cdots \ar[r]^{\phi_{m-1}} & 
    Y_m \ar[r]^{\sigma} \ar[d]^{h_m} & S \ar[d]^{h}\\
    X \ar[r]_{\psi_0} & Z_1 \ar[r]_{\psi_1} & Z_2 \ar[r]_{\psi_2} & 
    \cdots \ar[r]_{\psi_{m-1}} & Z_m \ar[r]_{\tau} & T }
  $$
  where the $\sigma$ and $\tau$ are contractions of fiber type with
  $\rho(Y_m/S) = \rho(Z_m/T) = 1$ and $h$ is finite with $\deg h =
  \deg f$.
\end{lem}

\begin{proof}
  The existence of the diagram follows again from
  Proposition~\ref{prop:etalebetweennonisom}. The fact that the
  general fiber of $\tau$ is Fano, hence simply connected, implies
  that $\deg h = \deg f$.
\end{proof}

In Sections~\ref{sec:MMPto1} and \ref{sec:MMPto2} we will consider the
cases where $\dim S$ is $1$ or $2$ separately. It will turn out in
either case that $S \cong T$ and that $h$ is \'etale.

\subsubsection{Minimal models over curves}\label{sec:MMPto1}

\begin{prop}\label{prop:FanoFibrationOverCurve} 
  In the setup of Lemma~\ref{lemnot:fibercontr}, suppose that $\dim S
  = 1$.  Then $S \cong \Alb(X)$, $h$ is \'etale and the fibration
  $\sigma$ is locally trivial in the analytic topology. The $\phi_i$
  are blow-ups of elliptic curves which are multi-sections over $S$.
\end{prop}

\begin{proof} 
  Apply Proposition~\ref{prop:etalebetweennonisom} to the two leftmost
  columns of Diagram~\eqref{eq:arrofofcontrs} and obtain
  $$
  \xymatrix{
    Y_m \ar[r]^{h_{m,0}} \ar[d]_{\sigma} & Z_{m,1} \ar[r]^{h_{m,1}}
    \ar[d]^{\tau_1} & Z_{m,2} \ar[r]^{h_{m,2}} \ar[d]^{\tau_1} & \cdots \\
    S \ar[r]_{h'_0} & T_1 \ar[r]_{h'_1} & T_2  \ar[r]_{h'_2} & \cdots
  }
  $$
  where $T_1$ are curves and $\tau_1$ are contractions of fiber type.
  Again, all $h'_i$ are finite of degree $d$.
  
  As the general fiber of $\sigma$ is rationally connected, it is
  clear that $q(X) = q(Y_m) = g(S)$. Observe that $g(S) \geq 1$: if
  not, $S \cong \bP_1$, and the theorem of Graber-Harris-Starr
  \cite{GHS03} would imply that $X$ is rationally connected, hence
  simply connected, a contradiction. The same argumentation shows
  $g(T_i) \geq 1$ for all $i$. But then
  $$
  \deg K_S \geq d\cdot \deg K_{T_1} \geq \cdots \geq d^i\cdot \underbrace{\deg K_{T_i}}_{\geq 0}
  $$
  for all $i$, which is possible if and only if $S$ and all $T_i$
  are elliptic, and all $h'_i$ \'etale.
  
  Since the general fiber of $\sigma$ is Fano, it follows that the
  Albanese map $X \to \Alb(X)$ factors via $S$. Since fibers of $X\to
  S$ are connected and $S$ is already elliptic, $S$ is $\Alb(X)$. The
  same holds for any of the $T_i$.
  
  To show that $\sigma$ is locally trivial, choose an arbitrary point
  $s \in S$ and observe that for any $i$ and any point of $s'$ of
  $$
  F_i := (h'_i \circ \cdots \circ h'_0)^{-1}(h'_i \circ \cdots
  \circ h'_0(s)),
  $$
  the scheme-theoretic fibers $\sigma^{-1}(s)$ and
  $\sigma^{-1}(s')$ are isomorphic. Since the cardinality $\#F_i =
  d^i$ becomes arbitrarily large, this shows that all scheme-theoretic
  $\sigma$-fibers are isomorphic. In particular, $\sigma$ is smooth
  and thus, by \cite[Thm.~4.2]{Kodaira86}, locally trivial in the
  analytic topology.
  
  The description of $\phi$ is immediate from
  Proposition~\ref{prop:divisorialcontr}.
\end{proof}

\begin{rem} 
  If $\sigma: Y_m \to S$ is not a $\bP_2$-bundle, we can say a bit
  more about its structure.
  
  First suppose that the fibers $F$ of $\sigma$ are proper del Pezzo
  surfaces, i.e., $K_F^2 \leq 6$. Then we can apply
  \cite[Prop.~0.4]{PS98} and obtain a finite \'etale cover $\tilde S
  \to S$ such that the fiber product $\tilde Y_m = Y_m \times_S \tilde
  S$ contains a divisor $D$ such that $D \to \tilde S$ is a
  $\bP_1-$bundle and $D \cap \tilde F $ is a $(-1)$-curve for all
  fibers $\tilde F$ of $\tilde Y_m \to \tilde S$. Then we can consider
  the resulting space $\tilde Y_{m,1} \to \tilde S$ and repeat the
  process, ending up with a $\bP_2-$bundle over an \'etale cover of
  $S$.
  
  If $\sigma $ is $\bP_1 \times \bP_1-$bundle, then we find a finite
  \'etale cover $\tilde S \to S$ such that $\tilde Y_m$ has a
  $\bP_1-$bundle structure of a surface $W$ which is in turn again a
  $\bP_1-$bundle over $\tilde S$.
\end{rem}

\subsubsection{Minimal models over surfaces}\label{sec:MMPto2}

We are next studying the case $\dim S = 2$. As in case $\dim S = 1$,
the morphism $h$ will turn out to be étale.

\begin{lem}
  In the setup of Lemma~\ref{lemnot:fibercontr}, suppose that $\dim S
  = \dim T = 2$. Then the morphism $h$ is étale of degree $d$.
\end{lem}
\begin{proof}
  It suffices to observe that all fibers of $\tau$ are conics, hence
  simply connected.
\end{proof}

\begin{prop}\label{prop:FanoFibrationOverSurface} 
 In the setup of Lemma~\ref{lemnot:fibercontr}, suppose that $\dim S
  = \dim T = 2$ and that $\kappa (S) \geq 0$. Then $S \cong T$. A
  finite \'etale base-change makes $Y_m$ a conic bundle over an
  abelian surface or over a product $C \times E$ with $E$ elliptic
  curve, and $C$ a curve of genus at least two.  In particular, in the
  second case $S$ is an \'etale quotient of the product of an elliptic
  curve with a curve of general type. 

  The discriminant locus of the conic bundle $\sigma$ is either empty
  or a disjoint union of elliptic curves.
\end{prop}
\begin{proof}
  The morphisms $\sigma$ and $\tau$ are conic bundles by Mori theory.
  We observe that $S \cong T$, because both are mapped isomorphically
  onto the same irreducible component of the cycle space of $X$.  In
  fact, otherwise $X$ would have two different 2-dimensional families
  of rational curves. This would imply that $S$ is uniruled, hence
  $\kappa (S) = - \infty$.
  
  Since $h$ is not an isomorphism, \cite[Thm.~3.2]{Fu02} shows that
  there are two distinct cases:
  \begin{enumerate}
  \item $\kappa (S) = 0$ and $S$ is
    abelian or hyperelliptic;
  \item $\kappa (S) = 1$ and after finite \'etale cover, $ S = E
    \times C$ with $E$ elliptic and $g(C) \geq 2$.
  \end{enumerate}
  It remains to prove the last statement.  Let $\Delta_S\subset S$ be
  the discriminant locus of $\sigma$. Again the fact that $f$ is
  \'etale implies that $h^{-1}(\Delta_S)=\Delta_S$. Thus $\Delta_S$ is
  either empty or a disjoint union of elliptic curves.
\end{proof}

\begin{prop}\label{prop:FanoFibrationOverRuled}
  Suppose $\kappa (X) = - \infty$ and that $S$ and $T$ are surfaces
  with $\kappa (S) = \kappa (T) = - \infty$. Then $S$ and $T$ are
  ruled surfaces over an elliptic curve $B$ and one of the following
  two situations occurs.
  \begin{enumerate}
  \item\ilabel{il:4ConicBundle} The minimal model $Y_m$ is a proper
    conic bundle, and the discriminant locus consists of \'etale
    multi-sections of $S \to B$.  Supposing moreover that $X$ and all
    blow-downs $Y_i$ have only finitely many extremal rays, or that
    already $X$ is minimal, then, replacing $f$ by an iterate $f_k$ if
    necessary, we obtain $Y_m = Z_m$ and $S = T = \bP(\sF)$ with $\sF$
    a semi-stable rank 2-bundle over $B$.
  \item\ilabel{il:4ProjBundle} The minimal model $Y_m$ is of the form
    $Y_m = \bP(\sE)$ with $\sE$ a rank 2-vector bundle over $S$ with
    $c_1(\sE)^2 = 4c_2(\sE)$. In this case, one of the following two
    conditions is satisfied:
    \begin{enumerate}
    \item[(1)]\ilabel{il:4ProjBundleA} for all ample classes $H$, the
      bundle $\sE$ is not $H$-stable, or
    \item[(2)]\ilabel{il:4ProjBundleB} $Y_m = S \times_B \bP(\sE')$ with
      $\sE'$ a rank 2-bundle over $B$, $T = \bP(\sE')$ and $\sE =
      p^*(\sE')\otimes \sL$ where $\sL\in\mbox{\rm Pic}(S)$, and $p: S
      \to B$ denotes the projection.
    \end{enumerate}
  \end{enumerate}
\end{prop}

\begin{proof} {\em Step 1.} Since $S$ and $T$ have both negative
  Kodaira dimension and since none of them can be rational (otherwise
  $X$ would be rationally connected, hence simply connected), $S$ and
  $T$ are birationally ruled over curves $B$ respectively, $C$ of
  positive genus. Let $p: S \to B$ resp. $q: T \to C$ denote the
  projections.  Then we observe that
  $$
  q(X) = q(Y_m) = q(S) = g(B) \text{\quad and \quad} q(X) = q(Z_m) =
  q(T) = g(C).
  $$
  The \'etale morphism $h: S \to T $ from Lemma
  \ref{lemnot:fibercontr} immediately yields an \'etale morphism $h':
  B \to C$ of degree $d \geq 2$. From Riemann-Hurwitz we deduce that
  $B$ and $C$ must be elliptic. Hence the composed maps $\alpha: X \to
  B$ and $\beta: X \to C$ coincide both to the Albanese map of $X$,
  which shows that $B = C$. As in the proof of
  Proposition~\ref{prop:FanoFibrationOverCurve}, the \'etalit\'e of
  $h'$ then immediately implies that $p$ and $q$ are submersions.
  Consequently $S$ and $T$ are minimal, i.e., ruled surfaces over the
  elliptic curve $B$.

  \medskip
  
  \noindent{\em Step 2.} Consider the case of \iref{il:4ConicBundle}
  and suppose now that the discriminant locus $\Delta_S \ne
  \emptyset$. So does $\Delta_T \ne \emptyset$ and both are disjoint
  unions of elliptic curves, i.e., multi-sections of $p$ and $q$.
  Assume that $X$ and all $Y_i$ have only finitely many extremal rays.
  Then $\overline{NE}(X/B)$ has only finitely many extremal rays.
  Hence we may pass to an iterate $f_k$ such that $Y_1 = Z_1$. Now we
  argue with $h_1$ instead of $f$ and proceed inductively to conclude
  $Y_m = Z_m$. A last application of this argument yields $S = T$,
  and, by Lemma \ref{lem:RuledSurfaces} below, $S$ is defined by a
  semi-stable vector bundle. This shows Claim~\iref{il:4ConicBundle}.

  \medskip
  
  \noindent{\em Step3.} We now consider case \iref{il:4ProjBundle} and
  suppose that $\Delta_S = \Delta_T = \emptyset$. Thus both $\sigma:
  Y_m \to S$ and $\tau: Z_m \to T$ are $\bP_1-$bundles. We have to
  distinguish two cases:
  \begin{enumerate}
  \item[(A)] $Y_m$ carries only one $\bP_1-$bundle structure, so $S =
    T$, or
  \item[(B)] $Y_m$ carries two $\bP_1-$ bundle structures.
  \end{enumerate}
  In both cases we can write $Y_m = \bP(\sE)$, where $\sE$ is a rank
  2-bundle over $S$. An explicit computation of Chern classes, using
  that $(K_{Y_m})^3=0$ and $(K_S)^2 = 0$, implies
  \begin{equation}\label{eq:ZeroDiscriminant}
    c_1(\sE)^2 = 4c_2(\sE). 
  \end{equation}
  In the case (A), if $\sE$ is stable with respect to some ample $H$,
  then Equation~\eqref{eq:ZeroDiscriminant} and \cite{AB97} or
  \cite[Thm.~3.7]{Takemoto72} imply that $\sE = p^*(\sE')\otimes \sL$,
  where $\sE'$ is a rank-$2$ bundle over $B$, and $\sL$ is a line
  bundle on $S$.  In particular, any fiber of $Y_m \to B$ is
  isomorphic to $\bP_1\times \bP_1$ and $Y_m$ has two different
  $\bP_1$-bundle structures contrary to our assumption.
  
  In case (B), any fiber $F$ of $Y_m \to B$ has two rulings, hence $F
  \simeq \bP_1 \times \bP_1$. If $b \in B$ is any point and $S_b$ the
  associated fiber, we can therefore normalize $\sE$ such that
  $\sE|_{S_b} = \sO^2_{S_b}$. It is then possible to write $\sE =
  p^*(\sE')$ with a rank-$2$ bundle $\sE'$ over $B$. Hence
  $$
  Y_m = S \times_B \bP(\sE')
  $$
  and $T = \bP(\sE')$. The \'etale map $h: S \to T$ means that writing
  $S = \bP(\sG), $ we have $\sG = j^*(\sE')$ up to a twist with a line
  bundle. This shows \iref{il:4ProjBundle} and ends the proof.
\end{proof}

We prove next the technical result in dimension 2 which was used above.

\begin{lem}\label{lem:RuledSurfaces} 
  Let $F:=\mathbb{P}(\sF)\to E$ be a ruled surface over an elliptic
  curve. If $S$ has a non-trivial \'etale endomorphism, then $\sF$ is
  semi-stable.
\end{lem}

\begin{proof}
  We know from \cite{Nak02} that $F$ has a non-trivial endomorphism.
  Moreover, if $\sF$ is indecomposable, then it was proved again in
  \cite{Nak02} that $F$ has an \'etale endomorphism. We have thus to
  analyze the decomposable case. After normalization \cite{Ha77}, we
  can assume $\sF\cong\sO_E\oplus L$, with $\deg(L)\le 0$. Then $\sF$
  is semi-stable if and only if $\deg(L)=0$. If $S$ has an \'etale
  endomorphism, it must have degree one on the fibers, and thus it
  comes from an endomorphism $\psi$ of $E$. The compatibility
  condition is that $\psi^*\sF\cong \sF\otimes \sL$ for some $\sL\in
  \Pic(E)$. This is possible only if either $\sL\cong\sO_E$ and
  $\psi^*L\cong L$, or $\sL\cong L^{-1}$ and $\psi^*L\cong L^{-1}$. In
  both cases, we obtain $\deg(L)=0$. The decomposable case really
  occurs, as exemplified by the trivial $2$-bundle on $E$.
\end{proof}

\subsubsection{Proper conic bundles with \'etale endomorphisms.}

In this part we give non-trivial examples of proper conic 
bundles with \'etale endomorphisms. The rough construction
idea is the following: it is clear that there are conic bundles 
with large relative Picard number having endomorphisms. 
We start with such a conic bundle, then we
try the drop the second Betti number by factorization.

\begin{prop}
  Let $B$ be an arbitrary curve, and $E$ be an elliptic curve, and
  denote $Y=E\times B$. Let $k$ be an odd positive integer. Then there
  exists a smooth threefold $X$ and a morphism $\phi : X \to Y$ with
  $\rho(X/Y)=1$, which realizes $X$ as a proper conic bundle with
  reduced, but not always irreducible fibers. Further, there exists an
  \'etale endomorphism $f$ of $X$ of degree $k^2$ making the following
  diagram commute:
  $$
  \xymatrix{
    X \ar[r]^{\phi} \ar[d]_{f} & Y  \ar[d]^{h}\\
    X \ar[r]_{\phi} & Y }
  $$
  where $h$ denotes the multiplication by $k$ on the first factor.
\end{prop}

\begin{proof}
  \emph{Step 1:} As a first step, we construct a proper conic bundle
  $S\to B$ with reduced, but not always irreducible fibers that
  carries a $B$-involution which interchanges the components of the
  reducible fibers.

  Let $D$ be an effective, reduced divisor on the curve $B$ such that
  $\sO_B(D)$ is divisible by two in $\mbox{\rm Pic}(B)$.  Consider
  $L\in \mbox{\rm Pic}(B)$ with $L^{\otimes 2}=\sO_B(D)$ and $s\in
  H^0(B,\sO_B(D))$ such that $s$ vanishes precisely along $D$. Denote
  $\mathbb{L}$ the total bundle space of $L$. It is known that the
  subvariety
  $$
  C:=\{x\in \mathbb{L}, x^{\otimes 2}=s\}
  $$
  is a smooth double covering of $B$ under the restriction of the
  bundle projection $\mathbb{L}\stackrel{\pi}{\to}B$.  The
  intersection between the inverse image of $D$ under $\pi$ and $C$
  coincides set-theoretically with the ramification divisor $D'$ of
  the covering $C \to B$.
  
  We compactify $\mathbb{L}$ to the ruled surface over $B$
  $$
  S:=\mathbb{P}(\sO_B\oplus L),
  $$
  where the projectivization is taken in the usual geometric sense,
  opposite to Grothendieck's. We keep the notation $\pi$ for the
  projection $S\to B$.  The natural morphism between vector bundles
  $$
  \left[
    \begin{array}{cc}
      0 & 1\\
      s & 0
    \end{array}
  \right] :\sO_B\oplus L\to L\oplus L^{\otimes 2}
  $$
  yields to a rational involution $S\dashrightarrow S$ which is
  relative over $B$ and whose indeterminacy locus is the ramification
  divisor $D'$ of $C\to B$. We consider the blow-up $\widetilde{S}$ of
  $S$ in the points of $D'$ so that the rational involution of $S$
  lifts to a morphism $\widetilde{S}\to S$. Note that the strict
  transforms of the fibers through the points of $D'$ will be
  contracted by this morphism, a fact which eventually proves that the
  original rational involution factors through a regular involution
  $\iota:\widetilde{S}\to \widetilde{S}$. By construction, the fixed
  locus of the involution $\iota$ coincides with the strict transform
  of $C$, and $\iota$ exchanges the irreducible components of each of
  the inverse images of fibers through the points of $D'$.

  \medskip
  
  \emph{Step 2:} Now $X$ is constructed as a suitable quotient of
  $S\times E$. We write the elliptic curve $E$ as
  $E=\mathbb{C}/\Lambda$, with $\Lambda =\mathbb{Z} +
  \mathbb{Z}.\omega$, and we identify the involution $\iota$ with a
  $2$-torsion element of $E$, say $1/2$. Via the canonical projection
  $E[2]\to \langle \iota\rangle$, $1/2\mapsto \iota$, $\omega/2\mapsto
  id$, the group $E[2]\cong\mathbb{Z}_2\oplus\mathbb{Z}_2$ of all
  $2$-torsion elements of $E$ acts diagonally on the product
  $E\times\widetilde{S}$.  Since the action of $E[2]$ on $E$ by
  translations is with trivial isotropy groups, so is the action of
  $E[2]$ on the product $E\times\widetilde{S}$, hence the quotient
  $$
  X:=(E\times\widetilde{S})/E[2]
  $$
  is a smooth variety.

  Consider next the group of $2k$-torsion elements $E[2k]$, and note
  that, since $k$ is odd, the inclusion $E[2]\subset E[2k]$ gives a
  decomposition $E[2k]\cong E[2]\times E[k]$.  Similarly to above, we
  identify $\iota$ with $1/2\in E[2k]$, and consider the diagonal
  action of $E[2k]$ on the product $E\times\widetilde{S}$. The action
  of the component $E[k]$ on $\widetilde{S}$ is trivial, and we
  observe that
  $$
  (E\times\widetilde{S})/E[2k]\cong X. 
  $$
  The inclusion $E[2]\subset E[2k]$ yields to a $(k^2:1)$ \'etale
  covering
  $$
  (E\times\widetilde{S})/E[2k]\to(E\times\widetilde{S})/E[2]
  $$
  that descends to the natural $(k^2:1)$ covering
  $$
  E\times B\cong E/E[2k]\times B\to E/E[2]\times B\cong E\times B.
  $$
  which is given by multiplication with $k$. It is elementary to check
  that all components of $\phi$-fibers have linearly dependent
  homology classes. The assertion that $\rho(X/Y)=1$ follows.
\end{proof}

\part*{The vector bundle associated with an endomorphism}
\section{Positivity of vector bundles associated to Galois coverings}

This section is devoted to the study of positivity properties of
vector bundles coming from Galois covers.  The Galois condition will
be used in the following form.

\begin{lem}\label{lem:Galois}
  Let $f:Y\to X$ be a finite morphism of degree $d$ between
  irreducible reduced complex spaces. Then $f$ is Galois if and only
  if the fibered product $Z := Y\times_XY$ decomposes as follows
  \begin{equation}
    \label{eq:galois-decomposition}
    Y\times_XY=\mathop\bigcup\limits_{j=1}^dZ_j,
  \end{equation}
  where the restriction of the first projection to any of the $Z_j$ is
  biholomorphic to $Y$.
\end{lem}

\begin{proof}
  The proof is rather straightforward, and very likely well-known. We
  notice first that the image of the diagonal map $Y\to Y\times_XY$ is
  one component $Z_1$ of $Y\times_XY$ which projects isomorphically to
  $Y$. Any automorphism $\sigma_i\in G(Y/X)$ naturally induces an
  automorphism in $\widetilde{\sigma}_j \in G\bigl( (Y\times_XY)/Y
  \bigr)$, acting on the second factor. The image
  $Z_j:=\widetilde{\sigma}_j(Z_1)$ is then an irreducible component of
  $Y\times_XY$ which projects biholomorphically to $Y$.
  
  If $f$ is Galois of degree $d$, then we get precisely $d$ such
  components. Since the degree of the first projection also equals
  $d$, we obtain a decomposition as in
  \eqref{eq:galois-decomposition}.

  Conversely, if $Z$ decomposes as in \eqref{eq:galois-decomposition},
  we can define elements in $G(Y/X)$ using that all components $Z_j$
  are isomorphic.
\end{proof}

\begin{notation}\label{not:Lazarsfeld-sheaf}
  Given any finite flat morphism $f: Y \to X$ with smooth target $X$,
  we consider the vector bundle
  $$
  \sE_f := \left(\factor f_*(\sO_Y).\sO_X. \right)^*.
  $$
  Recall that the trace map gives a splitting $f_*(\sO_Y) \cong
  \sO_X \oplus \sE_f^*$.
\end{notation}

To define $\sE_f$ it is a priori it is not necessary to make any
assumption on the smoothness of $Y$, nor of its components. Moreover,
flatness is preserved in some cases when components of $Y$ are
removed.

\begin{lem}\label{lem:flatness}
  Let $Z$ be a reduced projective variety and $f: Z \to Y$ a finite
  flat morphism, with smooth irreducible projective target $Y$. If $Z
  = \bigcup_{j=1}^d Z_j$ is the decomposition into irreducible
  components, and $Z' := \bigcup_{j=2}^d Z_j$, then the restriction
  $f' := f|_{Z'}$ is likewise flat.
\end{lem}

\begin{proof}
  We need to prove that $f'_*(\sO_{Z'})$ is locally free, or
  equivalently that the length $l(Z'_y)$ of all scheme-theoretic
  fibers $Z'_y := (f')^{-1}(y)$ are the same. We shall proceed by
  induction on $n = \dim Y$.
  
  If $n=1$, flatness is equivalent to $Y$ being dominated by any
  irreducible component of $Z'$, \cite[III~Prop. 9.7]{Ha77}. This
  condition is fulfilled.
  
  In general, take any point $y \in Y$ and consider a smooth connected
  hyperplane section $H \subset Y$ passing through $y$. The
  restriction $f|_H: Z_H \to H$ is obviously flat, and the restriction
  $f'|_H: Z'_H \to H$ is flat by induction hypothesis. Observe that
  $Z_H=\bigcup_jZ_{j,H}$, and that $Z_{j,H}$ surjects to $H$. Hence
  $$
  l(Z'_{H,x}) = \deg (f'_H) = \deg (f').
  $$
  Since $Z'_x = Z'_{H,x},$ we conclude.
\end{proof}

\begin{thm}\label{thm:SeveralComponents}
  Let $Z=\bigcup_{j=1}^dZ_j$ be a connected reduced projective variety
  of dimension at least 2, where $Z_j$ denote its irreducible
  components and $d\ge 2$.  Suppose that for all $k \not = j$, the
  intersection $Z_j \cap Z_k$ is either empty, or an ample Cartier
  divisor in $Z_j$.
  
  If $f:Z\to Y$ is a finite flat morphism to a smooth variety such
  that for any $j$, the restriction $f|_{Z_j}:Z_j\to Y$ is
  biholomorphic, then $\sE_f$ is ample.
\end{thm}

\begin{proof}
  We proceed by induction on the degree $d=\deg(f)\ge 2$.
  
  After renumeration of the irreducible components $Z_j$, we can
  assume without loss of generality that $Z'=\bigcup_{j\ge 2}Z_j$ is
  connected. Write $Z=Z_1\cup Z'$, and let $\ramification \subset Z_1$
  the scheme-theoretic intersection, $\ramification := Z_1 \cap Z'$.
  Since $\ramification$ is a union of ample divisors, it is ample.
  
  From Lemma~\ref{lem:flatness}, it follows that the induced map from
  $f':Z'\to Y$ is again flat. The Mayer-Vietoris sequence of the
  decomposition then reads as follows
  $$
  0\to \sO_Z\to \sO_{Z_1}\oplus \sO_{Z'}\to \sO_\ramification\to 0.
  $$
  Taking $f_*$, taking quotients by $\sO_Y$ and using that $Z_1$
  projects biholomorphically to $Y$, we obtain
  $$
  0\to \sE_f^*\to \sO_Y\oplus \sE_{f'}^*
  \stackrel{\beta}{\to}\sO_\ramification\to 0.
  $$
  Remark that $\beta|_{\sO_Y\oplus \{0\}}: \sO_Y \to \sO_\ramification$ is the
  obvious restriction map, so that
  $$
  \ker\left(\beta|_{\sO_Y\oplus \{0\}}\right)\cong \sO_Y(-\ramification).
  $$
  The Snake Lemma then yields a short exact sequence
  \begin{equation}\label{eq:galoisinducexten}
    0\to \sO_Y(-\ramification)\to \sE_f^* \to \sE_{f'}^* \to 0.    
  \end{equation}
  
  If $d=2$, then $\sE_{f'} = 0$, and
  Sequence~\eqref{eq:galoisinducexten} gives an isomorphism
  $\sE_f\cong \sO_Y(\ramification)$. Since the latter is ample, This
  settles the first induction step.
  
  If $d\ge 3$, the morphism $f'$ obviously satisfies all assumptions
  of Theorem~\ref{thm:SeveralComponents}, and $\sE_{f'}$ is ample by
  induction hypothesis. The dual of
  Sequence~\eqref{eq:galoisinducexten} then represents $\sE_f$ as an
  extension of ample bundles. It is thus ample, too.
\end{proof}

An almost identical argument, applied under weaker assumptions gives
rise to a nefness criterion.

\begin{prop}\label{prop:SeveralComponentsNef}
  Let $Z=\bigcup_{j=1}^dZ_j$ be a reduced projective variety of
  dimension at least 2, where $Z_j$ denote its irreducible components
  and $d\ge 2$.  Suppose that for all $k \not = j$, the intersection
  $Z_j \cap Z_k$ is either empty, or a nef Cartier divisor in $Z_j$.
  
  If $f: Z \to Y$ is a finite flat morphism to a smooth variety such
  that for any $j$, the restriction $f_j:=F|_{Z_j}:Z_j\to Y$ is
  biholomorphic, then $\sE_f$ is nef. \qed
\end{prop}

In view of Lemma~\ref{lem:Galois}, Theorem~\ref{thm:SeveralComponents}
applies to fibered products. The setup is the following. Suppose
$f:Y\to X$ is a finite Galois covering (of degree at least two) of
projective manifolds. Consider $Z=Y\times_XY$, and denote $g:Z\to Y$
the first projection. Then $\sE_f$ is ample if and only if $f^*\sE_f$
is ample. An obvious base-change formula gives $f^*\sE_f = \sE_g$.
Lemma~\ref{lem:Galois} applies, and we obtain a decomposition
$Z=\bigcup_{j=1}^dZ_j$, with all $Z_j$ isomorphic to $Y$.
Connectivity of $Z$ is studied in the following.

\begin{lem}\label{lem:Connectivity}
  With the notation of Lemma~\ref{lem:Galois}, the variety $Z$ is
  connected if and only if $f$ does not factor through $Y \to
  \widetilde{Y} \stackrel{\tilde{f}}{\to} X$, with
  $\widetilde{f}:\widetilde{Y}\to X$ is \'etale of degree at least
  two.
\end{lem}

\begin{proof} 
  Suppose first that $f$ factors through an \'etale covering as above.
  Then \cite[p. 63,~thm~17.4.1]{EGAIV4} shows that
  $\widetilde{Y}\times_X\widetilde{Y}$ is not connected, as the
  diagonal is one of its connected components. Since $Z$ obviously
  dominates $\widetilde{Y}\times_X\widetilde{Y}$, it cannot be
  connected.
  
  Conversely, suppose that $Z$ is not connected. Decompose $Z$ into
  its connected components, $Z =
  \widetilde{Z}_1\cup\dots\cup\widetilde{Z}_p$. We can assume that the
  diagonal $\Delta$ is contained in $\widetilde{Z}_1$. Denote
  furthermore $ \widetilde{Z}' = \widetilde{Z}_2 \cup \dots \cup
  \widetilde{Z}_p$.
  
  As in the proof of Lemma~\ref{lem:Galois}, we will identify an
  element in the Galois group $G(Y/X)$ with its induced automorphism
  of $Z$ without further mention. Consider the stabilizer of $\tilde
  Z_1$ in the Galois group,
  $$
  H:=\{\sigma\in G(Y/X) \,\,|\,\,
  \sigma(\widetilde{Z}_1)=\widetilde{Z}_1\}.
  $$
  The action of $H$ on $Z$ gives rise to an action of $H$ on
  $\widetilde{Z}_1$.  It is moreover clear that
  $\sigma(\widetilde{Z}')=\widetilde{Z}'$, for any $\sigma\in H$, so
  that $H$ acts also on $\widetilde{Z}'$.
  
  Set $\widetilde{Y} := Y/H$. Since $Y$ is irreducible, we obtain that
  $\widetilde{Y}$ is irreducible, too. We claim that the induced
  morphism $\widetilde{f} : \widetilde{Y} \to X$ is \'etale, in
  particular, we claim that $\widetilde{Y}$ is smooth. The proof is
  finished if this claim is shown.
  
  Following \cite[Thm.~17.4.1]{EGAIV4} again, to prove the claim, it
  suffices to prove that the diagonal $\Delta_{\tilde{Y}} \subset
  \widetilde{Y}\times_X\widetilde{Y}$ is a connected component of
  $\widetilde{Y}\times_X\widetilde{Y}$. To this end, remark that
  $$
  \widetilde{Y} \cong \factor \widetilde{Z}_1.H., \quad \text{and}
  \quad \widetilde{Y}\times_X\widetilde{Y} \cong \factor Z.H. \cong
  \factor \widetilde{Z}_1.H. \cup \factor \widetilde{Z}'.H..
  $$
  
  Since $\widetilde{Z}_1$ and $\widetilde{Z}'$ are disjoint and open,
  we obtain that $\widetilde{Z}_1/H$ and $\widetilde{Z}'/H$ disjoint
  and open in the quotient topology. Since $\Delta_{\widetilde{Y}} =
  \widetilde{Z}_1/H,$ we obtain the claim.
\end{proof}

We arrive at the main results of this section.

\begin{thm}\label{thm:ampleness1} 
  Let $f:Y\to X$ be a Galois covering of smooth varieties which does
  not factor through an \'etale covering of $X$, such that all
  irreducible components of the ramification divisor $\ramification$
  are ample on $Y$.  Then the bundle $\sE_f$ is ample.
\end{thm}

\begin{proof}
  Recall from Lemma~\ref{lem:Connectivity} that $Y\times_XY$ is
  connected. Now apply Theorem~\ref{thm:SeveralComponents}.
\end{proof}

\begin{cor}\label{cor:ampleness2} 
  Let $f:Y\to X$ be a Galois covering of degree at least two of smooth
  projective varieties with $\pi_1(X)=0$, and $\rho(Y)=1$. Then the
  bundle $\sE_f$ is ample. \qed
\end{cor}

The following example shows that the Galois condition in
Corollary~\ref{cor:ampleness2} is really necessary.

\begin{ex}
  In \cite[Example~2.1]{PS04} an example of a triple covering $f:Y \to
  X$ of Fano threefolds with $\rho = 1$ is established with the
  property that $\sE_f$ is not ample. Hence $f$ cannot be Galois;
  moreover the Galois group $G(Y/X)$ must be trivial, since $f$ cannot
  factor.
\end{ex}

\begin{prop}\label{prop:nefness} 
  Let $f:Y\to X$ be a Galois covering of smooth varieties, such that
  all the irreducible components of the ramification divisor
  $\ramification$ are nef on $Y$ and such that $\ramification$ is
  reduced. Then the bundle $\sE_f$ is nef.  \qed
\end{prop}

\section{Branched endomorphisms}
\label{sec:branched} \setcounter{equation}{0}

In this section, we study endomorphisms $f: X \to X$ with non-empty
branch locus. By Corollary~\ref{cor:3branchedgiveuniruled}, we know
that $X$ is uniruled. If $\rho(X) > 1$, one has to study the effect of
$f$ onto a Mori contraction. Here, however, we consider the case where
$\rho(X) = 1$, and $X$ is therefore Fano. We wish to address the
following well-known problem, at least under suitable assumptions.

\begin{problem}\label{prop:charcPn}
  Is $\bP_n$ the only Fano manifold with $\rho(X) = 1$ admitting an
  endomorphism of degree $\deg f > 1$?
\end{problem}
There are a number of cases where Problem~\ref{prop:charcPn} has a
positive answer.
\begin{itemize}
\item $X$ is 3-dimensional: \cite{ARV99, Schuhmann99,
    Hwang-Mok03}.
\item $X$ is rational homogeneous: \cite{PS89, HM99}.
\item $X$ is toric: \cite{OW02}.
\item $X$ contains a rational curve with trivial normal
  bundle: \cite[Cor.~3]{Hwang-Mok03}.
\end{itemize}
In this section, we give a positive answer in the following cases.
\begin{itemize}
\item $f$ is Galois and the branch locus $\branch$ is not stabilized
  by the action of $\Aut^0(X)$: Corollary~\ref{cor:611} and
  Proposition~\ref{prof:nodefandstab}.
\item $f$ is Galois, $X$ is almost homogeneous and $h^0(X, T_X) > n$:
  Corollaries~\ref{cor:611} and \ref{cor:asmall},
  Remark~\ref{cor:asmallx}.
\item $f$ is Galois and $X$ has a holomorphic vector field whose zero
  locus is not is not contained in the union of the branch loci of all
  the iterates of $f$: Corollary~\ref{cor:611} and
  Theorem~\ref{thm:VectorField}.
\end{itemize}
Here the condition that $f$ is Galois is only used to ensure that the
bundle $\sE_f$ is ample. In other cases we can guarantee that any
endomorphism $f: X \to X$ must have degree one.
\begin{itemize}
\item $X$ is of index $r \le 2$, and has a line $\ell$ not contained
  in the branch locus $R$ of $f$: Theorem~\ref{thm:smallIndex}.
\item $X$ satisfies the \emph{Cartan-Fubini Property} from
  Definition~\ref{def:CF}, $X$ is almost homogeneous and $h^0(X,
  T_X)>n$: Theorem~\ref{thm:CF}.
\end{itemize}
Numerous sub-cases, smaller results and variants are emphasized in the
text.

\subsection{Notation and assumptions}

In addition to the Assumptions~\ref{ass:endomorphism}, we fix the
following extra assumptions and notation throughout the present
Section~\ref{sec:branched}.
\begin{assumption}\label{ass:FanoPic1}
  Assume that $X$ is a Fano manifold of dimension $\dim X = n$, Picard
  number $\rho(X)=1$ and index $r$. We maintain the assumption that
  there exists an endomorphism $f: X \to X$ of degree $d \ge 2$.
\end{assumption}

\begin{notation}\label{not:muandother}
  Denote the ample generator of the Picard group by $\sO_X(1)$ and
  write
  $$
  f^*\bigl(\sO_X(1)\bigr) = \sO_X(\mu).
  $$
  We shall again consider the vector bundle
  $$
  \sE = \sE_f = \left(\factor f_*(\sO_X).\sO_X.\right)^*,
  $$
  that was already introduced in
  Notation~\ref{not:Lazarsfeld-sheaf}.  Recall that $\rank \sE =
  \deg(f)-1 = d-1$.
  
  If $k \in \mathbb N$ is any number, let $f_k = f\circ \cdots \circ
  f$ denote the $k$th iterate of $f$, and let $\sE_{f_k}$ be the
  associated vector bundle.
\end{notation}

\begin{rem}
  Equation~\eqref{eq:adjunction_for_f} immediately yields
  \begin{equation}\label{eq:BrachDiv}
    \sO_X(\ramification) = (1-\mu)\cdot K_X.
  \end{equation}
\end{rem}
\begin{rem}
  If $Z \subset X$ is a subvariety of pure dimension, define its
  degree as $\deg Z = Z.c_1\bigl(\sO_X(1)\bigr)^{\dim Z}$.  Comparing
  $c_1(\sO_X(1))$ with its pull-back, a standard computation of
  intersection numbers yields
  \begin{equation}
    \label{eq:degofpullback}
    \deg f^*(Z) = \mu^{\dim X - \dim Z} \cdot \deg Z.    
  \end{equation}
  In particular, if $Z$ is a point, we obtain $d = \mu^{\dim X}$.
\end{rem}

\subsection{Fano manifolds with small index}

We begin with the first case which is independent from the rest of the
Section. The proof uses methods from \cite{Schuhmann99}.

\begin{thm}\label{thm:smallIndex} 
  Let $X$ be Fano with $\rho (X) = 1$ and index $r \leq 2$. Let $f: X
  \to X$ be an endomorphism. Suppose that there is a line $\ell
  \subset X$ not contained in the branch locus $\branch$ of $f$. Then
  $\deg(f)=1$.
\end{thm}

\begin{proof} 
  We consider the local
  complete intersection curve $f^{-1}(\ell)$. As in
  \cite[p.~225]{Schuhmann99}, we have
  \begin{equation}\label{eq:f1}
    f^{-1}(\ell) \cdot \sO_X(1) = \frac{d}{\mu}    
  \end{equation}
  In particular, the curve $f^{-1}(\ell)$ contains at most
  $\frac{d}{\mu}$ irreducible components. If $\tilde \ell \to
  {f^{-1}(\ell)}$ is the normalization, this implies
  \begin{equation}\label{eq:f2}
    -2\frac{d}{\mu} \leq \deg \omega_{\tilde \ell}
  \end{equation}
  
  On the other hand, a standard computation, again taken from
  \cite{Schuhmann99}, shows that
  \begin{equation}\label{eq:Former2}
    \deg \omega_{f^{-1}(\ell)} = (r-2)d - r {\frac{d}{\mu}}.
  \end{equation}
  Recall that $ \deg \omega_{\tilde \ell} \leq \deg
  \omega_{f^{-1}(\ell)}$ with equality if and only if $f^{-1}(\ell)$
  is normal. Combining Equations~\eqref{eq:f2} and \eqref{eq:Former2},
  we obtain
  \begin{equation}\label{eq:f3}
    -2\frac{d}{\mu} \leq \deg \omega_{\tilde \ell} \leq \deg
    \omega_{f^{-1}(\ell)} = (r-2)d - r {\frac{d}{\mu}}.
  \end{equation}

  \subsubsection*{Case 1: $r = 1$}
  In this case Equation~\eqref{eq:f3} immediately gives $\mu = 1$.
  Hence $\deg f = 1$.

  \subsubsection*{Case 2: $r = 2$} 
  In this case, the left- and right hand side of \eqref{eq:f3} are
  equal. The equality of $-2\frac{d}{\mu} = \deg \omega_{\tilde \ell}$
  implies that $\tilde \ell$ contains precisely $\frac{d}{\mu}$
  components, each isomorphic to $\bP_1$. The equality in $\deg
  \omega_{\tilde \ell} \leq \deg \omega_{f^{-1}(\ell)}$ says that
  $f^{-1}(\ell)$ is smooth and \eqref{eq:f1} asserts that the
  components of $f^{-1}(\ell)$ are disjoint lines. To conclude, apply
  \cite[Lemma 4.2]{HM01}.
\end{proof}

\begin{rem}
  If $r = 2$ and there exists a line $\ell$ with trivial normal
  bundle, then \cite{Hwang-Mok03} proves that $\deg f = 1$. However without assumption on the normal bundle it
  might a priori happen that the deformations of a line $\ell$ are
  all contained in a divisor.
\end{rem}

\subsection{Positivity of $\sE_f$}

The next criteria for $X$ being isomorphic to the projective
space will generally rely on positivity properties of the bundle
$\sE_f$, or of its restrictions to rational curves.

The following Theorem~\ref{thm:5critB} and its immediate
Corollary~\ref{cor:5critA} are general recipies for answering
Problem~\ref{prop:charcPn} positively. Recall that a rational curve $C
\subset X$ with normalization $\eta: \bP_1 \to C$ is called
\emph{standard}, if $\eta^*(T_X) \cong \sO_{\bP_1}(2) \oplus
\sO_{\bP_1}(1)^{\oplus a} \oplus \sO_{\bP_1}^{\oplus b}$ for some $a,b
\geq 0$.

\begin{thm}\label{thm:5critB}
  Under the Assumptions~\ref{ass:FanoPic1}, $X \cong \bP_n$ if and
  only if there exists a number $k$ such that both of the following
  two conditions hold.
  \begin{enumerate}
  \item \ilabel{il:5critB1} There exists a covering family $(C_t)_{t
      \in T}$ of curves with $T$ an irreducible component of the Chow
    scheme, such that $C_t$ is a standard rational curve for general
    $t$ and such that $\sE_{f_k}|_{C_t}$ is ample for all $t$, and
  \item \ilabel{il:5critB2} $h^0\bigl(X,\, f_k^*(T_X)\bigr) >
    h^0\bigl(X,\,T_X\bigr)$.
  \end{enumerate}
\end{thm}

\begin{proof}
  First suppose that $X \cong \bP_n$, and let $k \geq 1$ be any
  number.  Condition~\iref{il:5critB1} follows from
  \cite[Thm.~6.3.55]{Laz-Positivity2}. Condition~\iref{il:5critB2}
  results from lifting back the Euler sequence via $f_k$.
  
  Now suppose conversely that a number $k$ is given such
  that~\iref{il:5critB1} and~\iref{il:5critB2} hold.  Consider the
  splitting
  \begin{align*}
    H^0 \bigl(X,\, (f_k)^*(T_X)\bigr) & \cong H^0(X,\, (f_k)_*(f_k)^*(T_X)\bigr) \\
    & \cong H^0\bigl(X,\, T_X \otimes (f_k)_*(\sO_X)\bigr) \\
    & \cong H^0\bigl(X,\, T_X \oplus (T_X \otimes \sE_{f_k}^*)\bigr) \\
    & \cong H^0\bigl(X,\, T_X\bigr) \oplus H^0(X,\, T_X \otimes \sE_{f_k}^*\bigr).
  \end{align*}
  By property~\iref{il:5critB2}, $\Hom_X\bigl(\sE_{f_k},\, T_X\bigr)
  \ne 0$, and we obtain a non-trivial map $ \sE_{f_k} \to T_X$.
  Theorem~\ref{thm:charactPn} then implies that $X \cong \bP_n$.
\end{proof}

The following corollary is a slightly less technical reformulation of
Theorem~\ref{thm:5critB}.

\begin{cor}\label{cor:5critA}
  Under the Assumptions~\ref{ass:FanoPic1}, $X \cong \bP_n$ if and
  only if there exists a number $k$ such that both of the following
  two conditions hold.
  \begin{enumerate}
  \item \ilabel{il:5critA1} The vector bundle $\sE_{f_k}$ is ample,
    and
  \item \ilabel{il:5critA2} $h^0\bigl(X,\, f_k^*(T_X)\bigr) >
    h^0\bigl(X,\,T_X\bigr)$.
  \end{enumerate}
\end{cor}
\begin{proof}
  If \iref{il:5critA1} and \iref{il:5critA2} hold,
  Theorem~\ref{thm:5critB} applies. If $X \cong \bP_n$,
  \iref{il:5critA2} follows from Theorem~\ref{thm:5critB} and
  \iref{il:5critA1} from \cite[Thm.~6.3.55]{Laz-Positivity2}.
\end{proof}

Condition~\iref{il:5critB1} in Theorem~\ref{thm:5critB} has a the
following useful reformulation.

\begin{lem}\label{lem:amplitude-and-connectedness}
  Let $C \subset X$ be a rational curve not in the branch locus $\branch$.
  The vector bundle $\sE_f|_C$ is ample if and only if $f^{-1}(C)$ is
  connected.
\end{lem}

\begin{proof}
Let $\tilde C$ be the normalization with its canonical morphism
  $\eta: \tilde C \to X$, and consider the fibered product $\tilde X
  := X \times_X \tilde C$ with its base change diagram.  
  $$
  \xymatrix {\tilde X \ar[r] \ar[d]_{\tilde f} & X \ar[d]^{f} \\
    \tilde C \ar[r]_{\eta} & X.}
  $$
  Then $\sE_{\tilde f} = \eta^*(\sE_f)$. The fibered product
  $\tilde X$ is connected if and only if $f^{-1}(C)$ is. Since $\eta$
  is finite, $\sE_{\tilde f}$ is ample if and only if $\sE_f|_C$ is.
  
  If $\tilde X$ is connected, \cite[Thm.~1.3]{PS00} asserts that
  $\sE_{\tilde f}$ is ample.  If $\tilde X$ is \emph{not} connected,
  the vector space of locally constant function on $\tilde X$ will be
  at least two-dimensional. Accordingly, $\sE_{\tilde f}^* = {\tilde
    f}_*(\sO_{\tilde X}) / \sO_{\tilde C}$ has at least one global
  section, and its dual cannot be ample.
\end{proof}

Combining Theorem~\ref{thm:5critB}, Theorem~\ref{thm:ampleness1} and
Lemma~\ref{lem:amplitude-and-connectedness}, we obtain the following.

\begin{cor}\label{cor:611}
  Let $X$ be Fano with $\rho(X) = 1$. Let $f: X \to X$ be an
  endomorphism with $\deg f > 1$. Then $X \simeq \bP_n$ if the
  following two conditions hold.
  \begin{enumerate}
  \item $f$ is Galois or there exists a standard rational curve
    $C \not \subset \branch$ such that $f^{-1}(C_t)$ is connected for all
    deformations $C_t$ of $C$ and no component of any $C_t$ is
    contained in $\branch$.
  \item $h^0\bigl(X,
    \, f^*(T_X)\bigr) > h^0\bigl(X,\,T_X\bigr)$.
  \end{enumerate}
  \qed
\end{cor}

\subsection{Manifolds with many vector fields} 
\label{sec:mfwithTXsects}

We will now consider Condition~\iref{il:5critB2} in
Theorem~\ref{thm:5critB}, which means that there are infinitesimal
deformations of $f_k$ which do not come from automorphisms of $X$. We
describe vector fields and group actions on manifolds for which
$h^0(X, T_X) = h^0(X,\, f^*(T_X))$, and give criteria for
Condition~\iref{il:5critB2} to hold.

\begin{prop}\label{prof:nodefandstab}
  If $h^0(X, T_X) = h^0(X, f^*(T_X))$, then there exists a surjective
  morphism of Lie groups $\eta: \Aut^0(X)\to \Aut^0(X)$ such that $f$
  is equivariant with respect to the natural action $\iota :
  \Aut^0(X)\times X \to X$ upstairs and the action $\iota \circ (\eta
  \times id)$ downstairs. In particular, the action of $\Aut^0(X)$
  stabilizes both the ramification and the branch loci.
\end{prop}
\begin{proof}
  Given a vector field $\vec v \in H^0(X, T_X)$, the assumption
  $h^0(X, T_X) = h^0(X, f^*(T_X))$ implies that there is a (unique)
  vector field $\alpha(\vec v)$ such that
  \begin{equation}
    \label{eq:infintfiberpres}
    df(\vec v) = f^*(\alpha(\vec v)),
  \end{equation}
  where $df: T_X \to f^*(T_X)$ is the differential of $f$.
  Equation~\ref{eq:infintfiberpres} implies that the action (upstairs)
  of the 1-parameter group associated with $\vec v$ is fiber
  preserving. Since $\vec v$ is arbitrary, the set of fiber preserving
  automorphisms is open in $\Aut^0(X)$. Since it is also closed, all
  automorphisms in $\Aut^0(X)$ are fiber preserving. The existence of
  $\eta$ then follows, e.g., from \cite[Prop.~1 on p.~14]{H-Oe}. The
  surjectivity is immediate.
\end{proof}

\begin{subrem}\label{rem:introalpha}
  The linear isomorphism $\alpha: H^0(X, T_X) \to H^0(X, T_X)$ is in
  fact the Lie-algebra morphism associated with the group morphism
  $\eta$.
\end{subrem}

\subsubsection{Zero loci of vector fields}

We will show that the assumption $h^0(X, T_X) = h^0(X,\, f^*(T_X))$
guarantees that the zero locus of any vector field $X$ maps to the
branch locus $\branch \subset X$ downstairs ---at least for a
sufficiently high iterate of $f$.  Let $\branch_k$ denote the
branch locus of $f_k: X \to X$. The zero-locus of a vector field $\vec
v \in H^0(X,\,T_X)$ will be denoted $Z(\vec v)$.

\begin{thm}\label{thm:VectorField}
  Under the Assumptions~\ref{ass:FanoPic1}, suppose that
  $h^0\bigl(X,\,T_X\bigr) = h^0\bigl(X,\,f^*(T_X)\bigr)$ and let $\vec
  v \in H^0\bigl(X,\,T_X\bigr) \setminus \{0\}$ be a non-zero vector
  field. Then there exists a number $k$ such that $f_k\bigl(Z(\vec
  v)\bigr) \subset \branch_k$.
\end{thm}

\begin{proof}
  We argue by contradiction and assume that $f_k\bigl(Z(\vec v)\bigr)
  \not \subset \branch_k$ for all $k$.  For any number $k$,
  choose an irreducible component $Z_k \subset Z(\vec v)$ such that
  \begin{enumerate}
  \item\ilabel{il:dick} $Z'_k := f_k(Z_k) \not \subset
    \branch_k$, and
  \item\ilabel{il:doof} $\dim Z_k$ is maximal among all irreducible
    components that satisfy \iref{il:dick}.
  \end{enumerate}
  Let $df_k: T_X \to f_k^*(T_X)$ be the differential of $f_k$.  Note that
  if $z \in Z_k$ is a general point, then $d$ is injective everywhere
  along the fiber $f_k^{-1}f_k(z)$.  In particular, $df_k$ is generically
  injective along all components of the subvariety $Z''_k :=
  f_k^{-1}(Z'_k)$, which is reduced, possibly reducible but of pure
  dimension.
  
  A repeated application of the assumption that
  $h^0\bigl(X,\,T_X\bigr) = h^0\bigl(X,\,f^*(T_X)\bigr)$ shows the
  existence of a vector field $\vec w_k \in H^0\bigl(X,\,T_X\bigr)$
  such that $df_k(\vec v)= f_k^*(\vec w_k)$.  The vector field $\vec w_k$
  will then vanish along $Z'_k$. Since $df_k$ is generically injective
  along all components of $Z''_k$, the vector field $\vec v$ will
  vanish along those components, i.e.~$Z''_k \subset Z(\vec v)$.
  
  A repeated application of Formula~\ref{eq:degofpullback} yields that
  $$
  \deg Z''_k = \deg f_k^*(Z'_k) = \mu^{k\cdot(\dim X - \dim Z)}
  \cdot \deg Z'_k \geq \mu^{k\cdot(\dim X - \dim Z)}.
  $$
  In particular, if $K := \limsup(\dim Z_k)$, then we obtain a
  subsequence $Z''_{k_i} \subset Z''_k$ such that
  \begin{itemize}
  \item the $Z''_{k_i}$ are reduced, possibly reducible subvarieties of
    $Z$ of pure dimension $K$, and
  \item the degree of the $Z''_{k_i}$ is unbounded: $\lim (\deg Z''_k)
    = \infty$.
  \end{itemize}
  This is clearly impossible.
\end{proof}

For special vector fields, it is not necessary to map the zero locus
down via $f_k$. Thus, a better statement holds.

\begin{thm}
  Under the Assumptions~\ref{ass:FanoPic1}, suppose that
  $h^0\bigl(X,\,T_X\bigr) = h^0\bigl(X,\,f^*(T_X)\bigr) \not = 0$.
  Then there exists a number $k$ and a vector field $\vec v_0$ such
  that $Z(\vec v_0) \subset \branch_k$.
\end{thm}

\begin{proof}
  Again, let $\alpha$ be the linear isomorphism discussed in
  Remark~\ref{rem:introalpha}. Let $\vec v_0$ be an eigenvector of
  $\alpha$ and observe that $Z(\vec v_0) = Z(\alpha(\vec v_0)) =: Z$.
  
  By Theorem~\ref{thm:VectorField}, the proof is finished if we show
  that $Z = f(Z)$, thus $Z = f_k(Z)$ for any $k$. For that, it
  suffices to note that the equality $d(\vec v_0) = f^*\alpha(\vec
  v_0)$ implies that $f$ is equivariant with respect to the flow of
  the vector fields $\vec v_0$ (upstairs) and $\alpha(\vec v_0)$
  (downstairs).
\end{proof}

\subsubsection{Almost homogeneous manifolds}

For the following recall that a manifold $X$ is \emph{almost
  homogeneous} if the automorphism group acts with an open orbit. Let
$E \subset X$ denote the complement of this open orbit; we call $E$
the exceptional locus of the almost homogeneous manifold.
Equivalently, $E$ is the minimal set such that $T_X|_{X\setminus E}$
is generated by global holomorphic vector fields on $X$. If $X$ is
homogeneous, the set $E$ is empty.

In this setup, Proposition~\ref{prof:nodefandstab} has two immediate
corollaries.

\begin{cor}\label{cor:apostata}
  Under the Assumptions~\ref{ass:FanoPic1}, let $X$ be an almost
  homogeneous projective manifold. Let $E_d \subset E$ be the
  (reduced, possibly empty) codimension-1 part of $E$ and write $E_d =
  \sO_X(a)$. Let $r$ be the index of $X$ and $\mu$ the number
  introduced in Notation~\ref{not:muandother}. If $\mu(r-a) \geq r$,
  then $h^0\bigl(X,\,f^*(T_X)\bigr) > h^0\bigl(X,\,T_X\bigr)$.
\end{cor}
\begin{proof}
  We argue by contradiction and assume $h^0\bigl(X,\,f^*(T_X)\bigr) =
  h^0\bigl(X,\,T_X\bigr)$. Recall from Equation~\ref{eq:BrachDiv} on
  page~\pageref{eq:BrachDiv} that $c_1(\ramification) = (\mu-1)r$. By
  Proposition~\ref{prof:nodefandstab}, we have an inclusion of reduced
  divisors $\branch \subseteq E_d \subset X$. Since $\ramification$
  is a strict subdivisor of $f^*(\branch)$, we obtain the
  following.
  $$
  (\mu-1)r = c_1(\ramification) < c_1(f^*(\branch)) \leq
  c_1(f^*(E_d)) = \mu\cdot a \quad \Leftrightarrow \quad \mu(r-a) \geq r
  $$
\end{proof}

\begin{cor}\label{cor:asmall}
  In the setup of Corollary~\ref{cor:apostata}, if $r>a$, then
  $h^0\bigl(X,\,f_k^*(T_X)\bigr) > h^0\bigl(X,\,T_X\bigr)$ for $k \gg
  0$.
\end{cor}
\begin{proof}
  Choose $k$ large enough so that $\mu^k(r-a) \geq r$ and apply the
  argumentation of Corollary~\ref{cor:apostata} to the morphism $f_k$.
\end{proof}

\begin{subrem}\label{cor:asmallx}
  Note that $r>a$ holds automatically if $h^0(X,\,T_X) > \dim X$, for
  the following reason. If $h^0(X,\,T_X) > n$, then $E_d$ is contained
  in the intersection of two distinct anticanonical divisors, given by
  wedge products of the form $s_0 \wedge \ldots \wedge s_{n-1} = 0$ or
  $s_1 \wedge \ldots \wedge s_n = 0$. Thus $a < r$. 
  
  Notice also as a special case of ~\ref{thm:VectorField} that if $\ramification
  \not \subset E,$ then $h^0(f^*(T_X)) > h^0(T_X)$.
\end{subrem}

\subsection{Manifolds that satisfy a Cartan-Fubini condition}
\label{sec:Cartan-Fubini}

Let $X$ be a Fano manifold and $T \subset \RatCurves^n(X)$ be a
dominating family of rational curves of minimal degrees. If $x \in X$
is a general point, we can look at the set of distinguished tangent
directions,
$$
\mathcal C_x := \bigl\{ \vec v \in \mathbb P(T_X^*|_x) \,|\,
\exists \ell \in T \text{ such that } x \in \ell \text{ and } \vec v
\in \mathbb P(T_\ell|_x)\bigr\}.
$$
We refer to \cite{KS06, Hwa00, Kebekus02a} for details.

Hwang and Mok have shown that for many Fano manifolds of interest,
the set of distinguished tangent directions in an analytic neighborhood
of the general point gives very strong local invariants that determine
the manifold globally. More precisely, they have shown that many Fano
manifolds satisfy the following Cartan-Fubini property.

\begin{defn}[see \cite{HM01}]\label{def:CF}
  A $X$ Fano manifold with $\rho(X) = 1$ is said to satisfy the
  \emph{Cartan-Fubini property}, (CF) for short, if there exists a
  dominating family $T \subset \RatCurves^n(X)$ of rational curves of
  minimal degrees such that the following holds.
  
  If $X'$ is any other Fano manifold with $\rho(X') = 1$, $S \subset
  \RatCurves^n(X')$ any dominating family of rational curves of
  minimal degrees and $\varphi: U \to U'$ a biholomorphic map between
  analytic open sets that respects the varieties of minimal rational
  tangents associated with $T$ and $S$, then $\varphi$ extends to a
  biholomorphic map $\Phi: X \to X'$.
\end{defn}

We refer to \cite{HM01} for details and for examples of Fano manifolds
that satisfy (CF). It is conjectured that almost all Fano manifolds
will have that property. We mention two examples that will be of
interest to us.

\begin{prop}\label{prop:CFcriteria}
  Let $X$ be a projective Fano manifolds with $b_2(X)=1$, not
  isomorphic to the projective space. If one of the following holds,
  then the manifold $X$ satisfies (CF).
  \begin{enumerate}
  \item\ilabel{il:dick2} $\dim X \geq 3$, and $X$ is prime of index $>
    \frac{\dim X+1}{2}$, \cite[1.2, 1.5 and 2.5]{Hwa00}
  \item $X$ is rational homogeneous, \cite[(2) on
    p.~566]{HM01} \qed
  \end{enumerate}
\end{prop}

\begin{thm}\label{thm:CF}
  Under the Assumptions~\ref{ass:FanoPic1}, suppose that $X$ satisfies
  (CF) and that $X$ is almost homogeneous. Then $h^0(X,\,T_X) = n$ and
  both the branch-- and the ramification divisor of $f$ are contained
  in the exceptional locus $E \subset X$.
\end{thm}
\begin{proof}
  Since $X$ satisfies (CF), \cite[Cor.~1.5]{HM01} asserts that $\dim_f
  \Hom(X,X) = \dim \Aut^0(X)$.  Corollary~\ref{cor:asmall} and
  Remark~\ref{cor:asmallx} then show that $h^0(X,\,T_X) \leq n$. But
  since $X$ is almost homogeneous, we have equality.
  Proposition~\ref{prof:nodefandstab} tells us where the branch and
  ramification loci are.
\end{proof}

Since rational homogeneous manifolds $X$ always satisfy $h^0(X,T_X) >
\dim X $, we obtain the following result of Paranjap\'e and Srinivas as
an immediate corollary.

\begin{cor}[\protect{\cite{PS89, HM99}}]\label{cor:rathom}
  Under the Assumptions~\ref{ass:FanoPic1}, suppose that $X$ is
  rational-homogeneous. Then $X$ is isomorphic to the projective
  space. \qed
\end{cor}

\begin{thm}\label{thm:delPezzo}
  Let $X$ be a del Pezzo manifold of degree $5$ of dimension $n \geq
  3$. Then $X$ does not admit an endomorphism of degree $\deg f > 1$.
\end{thm}
\begin{proof} 
  If $\dim X = 3$, then the claim is shown in \cite{ARV99},
  \cite{Schuhmann99}, and in \cite{Hwang-Mok03}.
  
  If $\dim X > 3$, we claim that $X$ is almost homogeneous with
  $h^0\bigl(X,\, T_X\bigr) > \dim X$.
  Proposition~\ref{prop:CFcriteria}\iref{il:dick} and
  Theorem~\ref{thm:CF} then immediately show
  Theorem~\ref{thm:delPezzo}. To prove the claim, recall from
  \cite[Thm.~3.3.1]{IP99} that $\dim X \leq 6$. Moreover, if $\dim X =
  6$, then $X$ is the Grassmannian $G(2,5)$ and we are done by the
  preceding corollary. In the remaining cases $n = 4,5$, recall from
  \cite[Sects.~7.8--13]{Fujita81} that there exists a diagram of
  birational morphisms
  $$
  \xymatrix{& \tilde X \ar[dl]_{\sigma} \ar[dr]^{\phi} \\ X & &
    \bP_n}
  $$
  where $\phi$ is the blow-up of a subvariety $F$ which is entirely
  contained in a linear hypersurface $H \subset \bP_n$. Notice that
  the vector fields in $\bP_n$ fixing $H$ pointwise span $T_{\bP_n}$
  outside $H$, simply because $T_{\bP_n}(-1)$ is spanned. Since $F $
  is contained in $H,$ these vector fields extend to $\tilde X$ and
  span $T_{\tilde X}$ outside $\phi^{-1}(H)$. Moreover
  $$
  h^0(\tilde X, \, T_{\tilde X}) \geq h^0(\bP_n,\, T_{\bP_n}(-1)) = n+1.
  $$
  Since $\sigma$ has connected fibers, the vector fields on $\tilde
  X$ descent to $X$, so that $X$ is almost homogeneous with $h^0(X,\,
  T_X) > n$. This shows the claim and ends the proof of
  Theorem~\ref{thm:delPezzo}.
\end{proof}

\appendix

\section{Manifolds whose tangent bundles contain ample subsheaves}

The main result of this section is a slight generalization of a result
of Andreatta and Wisniewski \cite{AW01} that characterizes the
projective space as the only Fano manifold $X$ with $b_2(X)=1$ whose
tangent bundle contains an ample, locally free subsheaf.

\subsection{Setup and statement of result}

Throughout this section, we consider a setup where $X$ is a Fano
manifold and $T' \subset \RatCurves^n(X)$ a dominating family of
rational curves of minimal degrees ---again we refer to \cite{KS06,
  Hwa00, Kebekus02a} for details about this notion.

\begin{notation}
  There exists a natural morphism $\iota : T' \to \Chow(X)$. Let $T :=
  \overline{\iota(T')} \subset \Chow(X)$ be the closure of its image.
  If $t \in T$ is any point, let $C_t \subset X$ be the reduction of
  the associated curve. The curve $C_t$ will then be rational,
  reduced, but not necessarily irreducible.
\end{notation}

The following is the main result of this section. Its proof is given
in Section~\ref{sec:Aproof} below.

\begin{thm}\label{thm:charactPn}
  Let $X$ be a Fano manifold with with $\rho(X) = 1$ and $\sF\subset
  T_X$ be a subsheaf of rank $r$. Assume that there exists a
  dominating family of rational curves of minimal degrees, $T' \subset
  \RatCurves^n(X)$, such that for any point $t \in T \subset
  \Chow(X)$, the restriction $\sF\vert _{C_t}$ is ample.
  
  Then $X \simeq \bP_n$, and either $\sF \simeq \sO_{\bP_n}(1)^{\oplus
    r}$ or $\sF\simeq T_{\bP_n}$.
\end{thm}

\begin{subrem}
  Recall that a coherent sheaf $\sF$ is called ample if the line
  bundle $\sO_{\bP(\sF)}(1)$ is ample. Details concerning this notion
  are discussed, e.g., in \cite{Ancona82}.
\end{subrem}

\begin{cor} 
  Let $X$ be a projective manifold with $\rho(X) = 1$ and $\sF \subset
  T_X$ an ample subsheaf. Then $X \simeq \bP_n$. \qed
\end{cor}

\subsection{Proof of Theorem~\ref*{thm:charactPn}}
\label{sec:Aproof}

In the remainder of the present section, we will prove
Theorem~\ref*{thm:charactPn}. For the reader's convenience, we
subdivide the proof into a number of fairly independent steps.

\subsubsection*{Step 0: Setup of notation}

If $t\in T$ is any point, the associated curve $C_t$ is a union of
irreducible, rational curves. Let $f_t: \coprod \bP_1 \to X$ be the
normalization morphism. We will consider the determinant $\det(\sF) :=
\wedge^{[r]}\sF \subset \wedge^rT_X$.  Finally, set $n := \dim X$.

\begin{rem}
  The sheaf $\sF$ is a subsheaf of a torsion-free sheaf, and therefore
  torsion free itself. In particular, if $\Sing(\sF)$ is the singular
  locus of $\sF$, i.e., the locus where $\sF$ is not locally free,
  then $\codim_X \Sing(\sF) \geq 2$.
\end{rem}

To prove Theorem~\ref*{thm:charactPn}, we will need to use the
following description of curves associated with general points of $T$,
\cite[Thm.~1.2]{Hwa00}.

\begin{fact}
  A general point $t \in T$ corresponds to an irreducible rational
  curve $C$ with normalization $f_t: \bP_1 \to C$ such that
  \begin{equation}\label{eq:stdn}
    f^*(T_X) = \sO_{\mathbb P_1}(2) \oplus 
    \sO_{\mathbb P_1}(1)^p \oplus \sO_{\mathbb P_1}^{n-p-1}.
  \end{equation}
  A curve for which~\eqref{eq:stdn} holds, is called ``standard''.
\end{fact}

\subsubsection*{Step 1: The rank and the Chern class of $\sF$}

\begin{lem}\label{lem:rank-and-chern}
  Let $t\in T$ be a general point. Then either one of the following
  holds true.
  \begin{enumerate-c}
  \item\ilabel{caesar} $r=n$ and $c_1(\sF) \cdot C_t = n+1$, or
  \item\ilabel{kleopatra} $c_1(\sF) \cdot C_t = r$.
  \end{enumerate-c}
\end{lem}
\begin{proof}
  We argue by contradiction and assume that
  \begin{equation}\label{eq:slope}
    c_1(\sF) \cdot C_t \ne r \text{\quad and \quad} r<n.
  \end{equation}
  Now consider a general point $t \in T$ and the associated
  irreducible, reduced curve $C_t$. Since $\codim_X \Sing(\sF) \geq
  2$, recall from \cite[lem.~2.1]{Hwa00} that the curve $C_t$ does not
  meet the singular locus of $\sF$, so that $f_t^*(\sF)$ is a locally
  free sheaf on $\bP_1$.  The ampleness assumption and
  Equation~\ref{eq:stdn} then imply that
  $$ 
  f_t^*(\sF) = \sO_{\mathbb P_1}(2) \oplus \sO_{\mathbb P_1}(1)^{r-1}
  \text{\quad or \quad} \sO_{\mathbb P_1}(1)^r.
  $$
  In particular, we have that
  $$
  c_1(\sF) \cdot C_t = c_1\bigl(f^*_t(\sF)\bigr) \geq r.
  $$
  Equation~\eqref{eq:slope} then implies that $c_1(\sF) \cdot C_t >
  r$. The slope $\mu$ with respect to the class of $C_t$ therefore
  fulfills the inequality
  $$
  \mu(\sF) \geq \frac{r+1}{r}>\frac{n+1}{n}\ge \mu(T_X).
  $$
  In particular, $T_X$ is not semistable and $\sF$ is destabilizing.
  From \cite[Prop.~1 and Prop.~3]{Hwa98}, the maximal destabilizing
  subsheaf $\sH$ of $T_X$ must satisfy $\mu(\sH)\le 1$, contradicting
  $\mu(\sF)>1$. Hence $\sF$ has rank $n$ and $c_1(\sF) \cdot C_t =
  n+1$.
\end{proof}

\subsubsection*{Step 2: Proof in case~\iref{caesar}}

The following Proposition ends the proof of
Theorem~\ref{thm:charactPn} in case~\iref{caesar}.

\begin{prop}
  Let $t\in T$ be a general point. If $r=n$ and $c_1(\sF) \cdot C_t =
  n+1$, then $X \cong \bP_n$ and $\sF = T_X$.
\end{prop}
\begin{proof}
  The injection $\sF \to T_X$ immediately gives an injection
  $\det(\sF) \to \det(T_X)$. In particular, we have that $-K_X.C_t
  \geq n+1$ for all $t \in T$. In this setup, \cite{Kebekus02b} gives
  that $X \cong \bP_n$, that the family $T'$ is proper and that the
  associated curves $C_t$ are lines\footnote{The statement of
    \cite[Thm.~1.1]{Kebekus02b} assumes that the inequality
    $-K_X.\ell \geq n+1$ holds for any curve $\ell \subset X$.
    Observe, however, that the proof of \cite[Thm.~1.1]{Kebekus02b}
    only uses curves $\ell$ coming from a given dominating family of
    rational curves of minimal degrees.}.

  It remains to show that $\sF = T_{\bP_n}$. We argue by contradiction
  and assume that the inclusion sequence has non-zero cokernel $Q$,
  \begin{equation}\label{eq:fint}
    \xymatrix{0 \ar[r] & \sF \ar[r] & T_{\bP_n} \ar[r]^{\beta} & Q
      \ar[r] & 0.}
  \end{equation}
  The equality of Chern classes implies that the support $S$ of $Q$
  has codimension at least $2$.

  To start, we claim that for any point $x \in S$, $\rank Q_x = 1$ as
  an $\sO_{X,x}-$module. Indeed, take a general line $\ell$ passing
  through $x$ and restrict sequence~\eqref{eq:fint} to $\ell$:
  $$
  \xymatrix{\sF|_{\ell} \ar[r] & T_{\bP_n}|_{\ell}
    \ar[r]^(.4){\beta|_{\ell}} & Q|_{\ell} \ar[r] & 0.}
  $$
  If $Q_x$ had a rank larger than $1$, then the splitting
  $T_{\bP_n}|_{\ell} \cong \sO_{\bP_1}(2) \oplus \sO_{\bP_1}(1)^{n-1}$
  immediately implies that the sheaf $\sF_{\ell}$ cannot be ample.
  This shows the claim that $\rank Q_x = 1$.
  
  Now choose a point $x \in S$ and consider the surjective quotient
  map at $x$, $\beta|_x: T_{\bP_n}|_x \to Q|_x$.  Here $Q|_x$ is a
  $1-$dimensional vector space so that its kernel $\ker(\beta|_x)$ has
  dimension $n-1$. The set of lines
  $$
  M := \{ \ell \subset \bP_n \text{ a line }\,|\, x \in \ell \text{
    and } T_{\ell}|_x \subset \ker(\beta|_x) \}
  $$ 
  then has dimension $n-2$, and the associated lines cover a divisor
  in $\bP_n$. In particular, if $\ell \in M$ is a general element,
  then $\ell$ is not contained in $S$. We fix a general $\ell \in M$
  for the sequel.

  Let $Q'$ denote the trivial extension of $Q_x$ to $X$, let $\gamma:
  T_{\bP_n} \to Q'$ be the quotient map and let $\sG = \ker(\gamma)$
  be its kernel. Since $\ell \not \subset S$, we know that
  $\sF|_{\ell}/\tor$ is a subsheaf of $\sG|_{\ell}/\tor$. By choice of
  $\ell$, the composed map $T_{\ell} \to T_X|_{\ell} \to Q'$ vanishes,
  hence $T_\ell \cong \sO_{\bP_1}(2)$ is a subsheaf of
  $\sG|_{\ell}/\tor$.  But then $\sG|_{\ell}/\tor$ has a factor
  $\sO_{\bP_1}(2)$ and since $\sG|_{\ell}/\tor$ is a proper subsheaf
  of $T_{\bP_n}|_{\ell} \cong \sO_{\bP_1}(2) \oplus
  \sO_{\bP_1}(1)^{n-1}$ of rank $n$, it cannot be ample.  It follows
  that $\sF|_{\ell}$ cannot be ample, a contradiction. Thus $Q = 0$
  and $\sF = T_{\bP_n}$.
\end{proof}

\subsubsection*{Step 3: Proof in case~\iref{kleopatra}}

We compare the sheaf $\sF_1 := \factor \wedge^r \sF.\tor.$ with the
locally free sheaf $\det \sF = \wedge^{[r]} \sF$. Since $\sF_1$ is
torsion free and has rank one, recall from \cite[Lem.~1.1.8 on
p.~147]{OSS} there exists a subscheme $Z\subset X$ such that
\begin{equation}\label{eq:incup}
  \sF_1 \cong \sI_Z \otimes \det \sF.  
\end{equation}

\begin{lem}\label{lem:outsideZ}
  Let $t \in T$ be any point and consider the associated cycle $C_t$.
  If $C_t \not \subset Z$, then the cycle associated with $t$ is
  irreducible, reduced and $C_t \cap Z = \emptyset$.

  In particular, if $t \in T$ is any point, then either $C_t \subset
  Z$, or $C_t \cap Z = \emptyset$.
\end{lem}
\begin{proof}
  To show that the cycle associated with $t$ is irreducible and
  reduced, we argue by contradiction and assume that there exists a
  component $C'_t \subset C_t$ such that $C'_t \not \subset Z$ and
  such that $c_1(\sF).C'_t < r$. If $\eta : {\mathbb P}_1 \to X$ is
  the normalization of $C'_t$, then we have a natural, non-trivial
  inclusion
  \begin{equation}\label{eq:incl}
    \sF' := \left(\factor \eta^*(\sF_1).\tor.\right) \to \eta^*(\det \sF).
  \end{equation}
  More precisely, the isomorphism~\eqref{eq:incup} then shows that
  \begin{equation}\label{eq:inclx}
    \sF' \cong \sI_{\eta^{-1}(Z)} \otimes \eta^*(\det \sF).
  \end{equation}
  The morphism~\eqref{eq:incl} shows that the locally free sheaf
  $\sF'$ on $\mathbb P_1$ has rank $r$ and degree $\deg \sF' \leq \deg
  \eta^*(\det \sF) < r$, so $\sF'$ cannot be ample. On the other hand,
  $\sF'$ is a quotient of a wedge power of the pull-back of an ample
  sheaf under a finite morphism, hence ample. This contradiction shows
  $C'_t = C_t$.

  We can now assume that $C_t$ is irreducible and that
  $C_t.c_1(\sF)=r$. Again, let $\eta: \mathbb P_1 \to X$ be the
  normalization, and consider the morphism~\eqref{eq:incl}. Again,
  $\sF'$ is ample, but this time $\deg \sF' \leq \deg \eta^*(\det \sF)
  = r$. Obviously, if $\sF'$ is ample, then $\deg \sF' \geq r$.
  Isomorphism~\eqref{eq:inclx} asserts that this is the case if and
  only if $\eta^{-1}(Z) = \emptyset$, as claimed.
\end{proof}

\begin{cor}\label{cor:zempty}
  The subscheme $Z \subset X$ is empty. In particular, the family $T$
  is non-split, i.e., if $t \in T$ is any point, then the associated
  curve $C_t$ is irreducible and satisfies $C_t.c_1(\sF)=r$.
\end{cor}
\begin{proof}
  By Lemma~\ref{lem:outsideZ}, both $Z$ and $X\setminus Z$ are unions
  of curves of the form $C_t$.  Fix a general point $x \in X \setminus
  Z$ and consider the compact variety $Z^{[1]}$ filled up by all
  curves $C_t$ through $x$. Then consider the compact variety
  $Z^{[2]}$ filled up by all curves $C_t$ meeting $Z^{[1]}$ and so on.
  By Lemma~\ref{lem:outsideZ}, all $Z^{[k]}$ are contained in $X
  \setminus Z$. Now the sequence $(Z^{[k]})_{k \in \mathbb N}$ must
  stabilize at some $k_0$.

  Assume that $Z \not = \emptyset$. Then $Z^{[k_0]}$ is a proper
  subvariety of $X$. In that case, the family $(C_t)$ is not
  connecting, and we obtain an almost holomorphic map $X \dasharrow W$
  to a normal variety $W$ of positive dimension. This contradicts
  $\rho (X) = 1$ and shows the claim.
\end{proof}

\begin{cor}
  The double dual $\sF^{**}$ is locally free.
\end{cor}
\begin{proof}
  Following \cite[Satz~1.1, Rossi's Theorem]{GrauertRiemenschneider},
  there exists a sequence of blowups $\pi:\widehat{X}\to X$ such that
  $\pi^*(\sF)/\tor$ is locally free.
  \begin{align*}
  \pi^*\det(\sF) & \cong \pi^*(\wedge^r \sF / \tor)  \\
  & \cong \pi^*(\wedge^r \sF)/\tor && \text{by \cite[Satz~1.3]{GrauertRiemenschneider}} \\
  & \cong \wedge^r(\pi^* \sF)/\tor  \\
  & \cong \wedge^r(\pi^* \sF /\tor)/\tor \\
  & \cong \wedge^r(\pi^* \sF /\tor) \cong  \det(\pi^* \sF /\tor)&& \text{since $\pi^* \sF /\tor$ is locally free }
  \end{align*}
  In particular, if $\widehat{X}_x$ is any $\pi$-fiber, then the
  restriction $\det(\pi^* \sF /\tor)|_{\hat{X}_x} \cong
  \sO_{\hat{X}_x}$ is trivial. Since $(\pi^* \sF /\tor)|_{\hat{X}_x}
  \subset \pi^*(T_X)|_{\hat{X}_x} \cong \sO_{\hat{X}_x}^{\oplus \dim
    X}$, it follows that $(\pi^* \sF /\tor)|_{\hat{X}_x} \cong
  \sO_{\hat{X}_x}^{\oplus r}$. 
  
 Consequently, there exists a vector bundle
  $\sF'$ on $X$ such that $\pi^* \sF /\tor \cong \pi^*(\sF')$.  The
  sheaves $\sF^{**}$ and $\sF'$ are isomorphic outside the singular
  locus of $\sF$, which of codimension at least two, hence $\sF^{**} =
  \sF'$ is locally free.
\end{proof}

\begin{cor}
  The sheaf $\sF$ is locally free.
\end{cor}
\begin{proof}
  Consider the natural exact sequence
  \begin{equation}\label{eq:Q}
    0 \to \sF \to \sF^{**} \to Q \to 0,    
  \end{equation}
  where $Q$ is a torsion sheaf. We need to show that $Q = 0$.
  
  Let $t \in T$ be any point and $C_t$ the associated curve which, by
  Corollary~\ref{cor:zempty}, is irreducible. We claim that either
  $C_t \subset \Supp(Q)$, or that $C_t$ is disjoint from $\Supp(Q)$.
  Indeed, if $C_t$ would meet $\Supp(Q)$ in a finite non-empty set,
  and if $\eta : \bP_1 \to X$ is the normalization morphism,
  Sequence~\eqref{eq:Q} pulls back to
  $$
  \eta^*(\sF)/\tor \to \underbrace{\eta^*(\sF^{**})}_{\cong
    \sO_{\bP_1}(1)^{\oplus r}} \to \underbrace{\eta^*(Q)}_{\not \cong
    0} \to 0.
  $$
 In particular, $\eta^*(\sF)/\tor$ is a
  locally free strict subsheaf of $\sO_{\bP_1}(1)^{\oplus r}$ and
  cannot be ample. This shows the claim.
  
  As in the proof of Corollary~\ref{cor:zempty}, consider the
  varieties $Z^{[k]}$ covered by curves that, observe that the
  $Z^{[k]}$ do not intersect $\Supp(Q)$ and conclude that $\Supp(Q)$
  is empty.
\end{proof}

\medskip

We have shown that $\sF$ is locally free. To finish the proof of
Theorem~\ref{thm:charactPn}, apply \cite{AW01}.

\end{document}